\documentclass[reqno,11pt,pdftex]{amsart}
\usepackage[top=25truemm,bottom=25truemm,left=25truemm,right=25truemm]{geometry}
\usepackage{amsmath,amsthm,amssymb,amsfonts,amscd,mathtools,color,tikz}
\usepackage[all]{xy} 
\usepackage[colorlinks,linkcolor=black,citecolor=black,urlcolor=blue]{hyperref}
\numberwithin{equation}{section}
\theoremstyle{plain} 
	\newtheorem{thm}{Theorem}[section]
	\newtheorem*{thm*}{Theorem}
	\newtheorem{cor}[thm]{Corollary}
	\newtheorem{lem}[thm]{Lemma}
	
	\newtheorem{prop}[thm]{Proposition}
	\newtheorem{conj}[thm]{Conjecture}

\theoremstyle{definition}
	\newtheorem{defn}[thm]{Definition}

\theoremstyle{remark}
	\newtheorem{rem}[thm]{Remark}
	
	\newtheorem*{pf}{Proof}

\def\CC{{\mathbb C}}

\def\LL{{\mathbb L}}

\def\QQ{{\mathbb Q}}
\def\RR{{\mathbb R}}
\def\TT{{\mathbb T}}

\def\XX{{\mathbb X}}
\def\ZZ{{\mathbb Z}}

\def\B{{\mathcal B}}

\def\D{{\mathcal D}}
\def\E{{\mathcal E}}
\def\F{{\mathcal F}}

\def\I{{\mathcal I}}

\def\L{{\mathcal L}}

\def\O{{\mathcal O}}
\def\P{{\mathcal P}}

\def\S{{\mathcal S}}
\def\T{{\mathcal T}}

\def\X{{\mathcal X}}
\def\Z{{\mathcal Z}}
\def\Aut{{\rm Aut}}
\def\Br{{\rm Br}}

\def\FEC{{\rm FEC}}

\def\Hom{{\rm Hom}}

\def\ST{{\rm ST}}
\def\Stab{{\rm Stab}}
\def\rad#1{\mathrm{rad}(#1)}
\def\mod{{\rm mod}}

\def\h{{\mathfrak h}}

\def\p{{\partial}}
\def\br#1{\langle #1 \rangle}
\def\ns{{\nabla}\hspace{-1.4mm}\raisebox{0.3mm}{\text{\footnotesize{\bf /}}}}
\begin{document}
\title{Twist automorphism for a generalized root system of affine ADE type}
\date{\today}
\author{Takumi Otani}
\address{Yau Mathematical Sciences Center, Tsinghua University, Haidian District, Beijing, China}
\email{takumi@tsinghua.edu.cn}

\maketitle
\begin{abstract}
For a generalized root system of affine ADE type, we introduce a twist automorphism.
We prove that the Dubrovin--Zhang extended affine Weyl group is isomorphic to our (modified) extended affine Weyl group, which is an extension of the affine Weyl group by the twist automorphism.
We also show that the number of root bases with a Coxeter transformation modulo the twist automorphism is equal to the degree of the Lyashko--Looijenga map of the Frobenius manifold constructed by Dubrovin--Zhang. 
As analogues of the extended affine Weyl group, we define an extended Artin group and an extended Seidel--Thomas braid group.
We study the relationship between the extended affine Weyl group and the extended Seidel--Thomas braid group. 
\end{abstract}
\section{Introduction}
A {\em generalized root system} is defined as a pair consisting of a root system and a Coxeter transformation.
This notion was introduced by Kyoji Saito \cite{S1, S3} in his study of the Milnor lattice of an isolated hypersurface singularity.
A generalized root system is an important tool to study a triangulated category.
Furthermore, a {\em Frobenius manifold} associated with a generalized root system is conjecturally related to the space of stability conditions on a triangulated category (cf.~\cite{Br2, BQS, IQ, IQ2, IOST}).

Dubrovin \cite{D} constructed a Frobenius manifold for a (generalized) finite root system, which was originally found by \cite{SYS, S2} as a {\em flat structure}.
For an affine root system (with an additional choice of a vertex of the affine Dynkin diagram), Dubrovin--Zhang \cite{DZ} introduced an {\em extended affine Weyl group} and constructed a Frobenius manifold via the invariant theory of the group.
In this paper, we call the group the {\em Dubrovin--Zhang extended Weyl group}.
Based on the correspondence between root systems of affine $A$ type and the genus $0$ Hurwitz moduli spaces, their extended affine Weyl group was introduced as the monodromy group of the Frobenius manifold.

In this paper, for a generalized root system of affine ADE type, we introduce a twist automorphism and give an intrinsic definition of an extended affine Weyl group.
In order to define the twist automorphism and the extended affine Weyl group, it is necessary to consider a generalized root system.
\bigskip
Let $A = (a_1, a_2, a_3)$ be a tuple of positive integers satisfying $1/a_1 + 1/ a_2 + 1/ a_3 > 1$ and set $\mu_A = a_1 + a_2 + a_3 - 1$.
We consider a generalized root system $(\widetilde{R}_A, \widetilde{\mathbf{c}}_A) = (\widetilde{L}, \widetilde{I}, \widetilde{\Delta}_\mathrm{re}, \widetilde{\mathbf{c}}_A)$ of affine ADE type associated with $A$.
One can define the notion of an Euler form for a generalized root system. 
In the case of $(\widetilde{R}_A, \widetilde{\mathbf{c}}_A)$, there exists a unique Euler form $\chi \colon \widetilde{L} \times \widetilde{L} \longrightarrow \ZZ$ (see Lemma \ref{lem : uniqueness of the Euler from}).
For a generator $\delta \in \widetilde{L}$ of the radical $\mathrm{rad}(\widetilde{I}) \cong \ZZ$, the {\em twist automorphism} $t_\delta \in \Aut_\ZZ (\widetilde{L}, \widetilde{I})$ for $(\widetilde{R}_A, \widetilde{\mathbf{c}}_A)$ along $\delta$ is define by 
\begin{equation*}
t_\delta (\lambda) \coloneqq \lambda - \chi(\delta, \lambda) \delta, \quad \lambda \in \widetilde{L}.
\end{equation*}
Then, we define the {\em extended affine Weyl group} $\widetilde{W}_A^\mathrm{ext}$ as the subgroup of $\Aut_\ZZ(\widetilde{L}, \widetilde{I})$ generated by the affine Weyl group $\widetilde{W}_A$ and $T_A \coloneqq \br{t_\delta}$ (Definition \ref{defn : extended affine Weyl group}).

Let $\widetilde{\h}_A \coloneqq \Hom_\ZZ(\widetilde{L}, \CC)$ and define subsets $X_A$ and $X_A^\mathrm{reg}$ by
\begin{equation*}
X_A \coloneqq \widetilde{\h}_A \setminus H_\delta, \quad X_A^\mathrm{reg} \coloneqq X_A \setminus \bigcup_{\alpha \in \widetilde{\Delta}_\mathrm{re}} H_\alpha,
\end{equation*}
where $H_\lambda \coloneqq \{ x \in \widetilde{\h}_A \mid x(\lambda) = 0 \}$. 
Denote by $\pi_X \colon \X_A \longrightarrow X_A$ the universal covering space of $X_A$ and put $\X_A^\mathrm{reg} \coloneqq \pi_X^{-1}(X_A^\mathrm{reg})$.
Then, the {\em modified extended affine Weyl group} $\widehat{W}_A$ is defined to be the fiber product of $\widetilde{W}_A^\mathrm{ext}$ and the diagonal subset $C_A \subset T_A \times \pi_1(X_A)$ over $T_A$.
The group $\widehat{W}_A$ is isomorphic to $\widetilde{W}_A^\mathrm{ext}$ as a group, however their natural actions on $\X_A$ are different.
On the other hand, the Dubrovin--Zhang extended affine Weyl group $\widehat{W}_A^\mathrm{DZ}$ is defined by actions on $\widehat{\h}_A$, where $\widehat{\h}_A \coloneqq \h_A \times \CC$ and $\h_A = \Hom_\ZZ (\widetilde{L} / \mathrm{rad}(\widetilde{I}), \CC)$ (see Definition \ref{defn : Dubrovin--Zhang extended affine Weyl group}).
The following theorem is the first main result:
\begin{thm}[Proposition \ref{prop : extended affine Weyl groups} and Theorem \ref{thm : main 1}]\label{thm : introduction 1}
There is a group isomorphism $\widehat{W}_A \cong \widehat{W}_A^\mathrm{DZ}$.
Moreover, there exists a $\widehat{W}_A$-equivariant biholomorphism $\widehat{\varphi} \colon \X_A \longrightarrow \widehat{\h}_A$.
\end{thm}
As a corollary, one can obtain the Frobenius structure constructed by Dubrovin--Zhang on the orbit space $\X_A / \widehat{W}_A$.
In this paper, we study some invariants of the Frobenius manifold from the view point of the generalized root system with the twist automorphism.
For the Frobenius manifold $\X_A / \widehat{W}_A$, we naturally consider the Lyashko--Looijenga map $LL \colon \X_A / \widehat{W}_A \longrightarrow \CC^{\mu_A}$, which is a ramified covering.
For the generalized affine root system $(\widetilde{R}_A, \widetilde{\mathbf{c}}_A)$, put 
\begin{equation*}
\B_A \coloneqq \{ (\alpha_1, \dots, \alpha_{\mu_A}) \mid \widetilde{\mathbf{c}}_A = r_{\alpha_1} \cdots r_{\alpha_{\mu_A}}, \ \text{$B = \{ \alpha_1, \dots, \alpha_{\mu_A} \}$ is a root basis of $\widetilde{R}_A$} \}.
\end{equation*}
Denote by $e (\B_A)$ the cardinality of $\B_A$ modulo $T_A$.
Otani--Shiraishi--Takahashi proved that the degree of the Lyashko--Looijenga map is equal to the number of full exceptional collections in the derived category of the extended Dynkin quiver of type $A$ modulo spherical twists (\cite[Theorem 1.2]{OST}).
Based on the result, we prove the following
\begin{thm}[Corollary \ref{cor : LL map and root bases}]\label{thm : introduction 2}
We have
\begin{equation*}
\deg LL = e ( \B_A ) = \frac{\mu_A!}{a_1! a_2! a_3! \chi_A}a_1^{a_1}a_2^{a_2}a_3^{a_3}.
\end{equation*}
\end{thm}

As an analogue of the extended affine Weyl group, we consider the {\em extended Artin group} $\widehat{G}_A$ (Definition \ref{defn : extended Artin group}).
One can see that the extended Artin group $\widehat{G}_A$ is isomorphic to the fundamental group $\pi_1(\X_A^\mathrm{reg} / \widehat{W}_A)$ (Proposition \ref{prop : fundamental group and extended Artin group}).
Moreover, there exists a surjective homomorphism $\widehat{G}_A \longrightarrow \widehat{W}_A$ such that the isomorphism $\widehat{G}_A \cong \pi_1(\X_A^\mathrm{reg} / \widehat{W}_A)$ is compatible with the monodromy group of the Frobenius manifold $\X_A / \widehat{W}_A$ (see Corollary \ref{cor : extended Artin group and extended Weyl group}).

Let $\CC \widetilde{\TT}_A$ be a $\CC$-algebra associated with a quiver with relation $(Q_{\widetilde{\TT}_A}, \I_{\widetilde{\TT}_A})$ (Definition \ref{defn : octopus algebra}).
The derived category $\D_A \coloneqq \D^b \mod(\CC \widetilde{\TT}_A)$ admits a $1$-spherical object, and we denote by $\ST_1(\D_A)$ the Seidel--Thomas braid group of $\D_A$.
On the other hand, one can consider an $\XX$-Calabi--Yau triangulated category $\D^\XX_A$ and (the principal part of) the Seidel--Thomas braid group $\ST^\circ(\D^\XX_A)$ of $\D^\XX_A$.
There is a functor $\L \colon \D_A \longrightarrow \D^\XX_A$, which is called the Lagrangian immersion.
Then, the {\em extended Seidel--Thomas braid group} $\widehat{\ST}(\D^\XX_A)$ of $\D^\XX_A$ (Definition \ref{defn : extended Seidel--Thomas braid group}) is defined to be 
\begin{equation*}
\widehat{\ST}(\D^\XX_A) \coloneqq \big\langle \ST^\circ(\D^\XX_A), \ \ST^\L_1(\D_A) \big\rangle,
\end{equation*}
where $\ST^\L_1(\D_A)$ is the image of $\ST_1(\D_A)$ by the Lagrangian immersion.

\begin{thm}[Theorem \ref{thm : extended Artin group and extended ST group}]
There exists a surjective homomorphism $\widehat{\iota}_\XX \colon \widehat{G}_A \longrightarrow \widehat{\ST}(\D^\XX_A)$.
\end{thm}
As a consequence, we obtain the following diagram:
\begin{equation*}
\xymatrix{
\pi_1(\X_A^\mathrm{reg} / \widehat{W}_A) \ar@{->>}[rr] \ar[d]_\rho & & \widehat{\ST}(\D^\XX_A) \ar[d]^{\rho_\XX} \\
\widehat{W}_A \ar@{=}[rr] & & \widehat{W}_A 
}
\end{equation*}
where $\rho$ is the map defining the monodromy of the Frobenius manifold $\X_A / \widehat{W}_A$ and $\rho_\XX$ is defined in \eqref{eq : extended Braid group and extended affine Weyl group}.
We expect the map $\widehat{\iota}_\XX \colon \widehat{G}_A \longrightarrow \widehat{\ST}(\D^\XX_A)$ is an isomorphism (Conjecture \ref{conj : Artin and braid}).

As an analogue of \cite[Section 1.2]{IQ}, it is natural to expect that the orbit space $\X_A / \widehat{W}_A$ is biholomorphic to the quotient space of stability conditions on $\D_A$ by $\ST_1(\D_A)$.
It is also natural to expect that the regular part $\X_A^\mathrm{reg} \big/ \widehat{W}_A$ is biholomorphic to the quotient space of $q$-stability conditions on $\D^\XX_A$ by $\widehat{\ST}(\D^\XX_A)$ (see Conjecture \ref{conj : stability condition}).
These conjectures are proved for $A = (1, p, q)$ and $\XX = [N]$ with $N \ge 3$ (Proposition \ref{prop: affine A case}).
\bigskip
\noindent
{\bf Acknowledgement.} 
I would like to thank Akishi Ikeda and Yu Qiu for careful reading of this manuscript and sharing their insights on $\XX$-Calabi--Yau categories and Seidel--Thomas braid groups.
I also thank Ping He for his valuable comments. 
This work is supported by Beijing Natural Science Foundation Grant number IS24008.
\section{Generalized roost system}\label{sec : Generalized roost system}
We recall the notion of a generalized root system in this section.
In order to define the twist automorphism, we need to consider a generalized root system.
We refer the reader to \cite{H, K} for root systems and \cite{S1, S3, T, STW} for generalized root systems.
Throughout this paper, we always assume that our root systems are simply-laced.
\subsection{Root system}
Let $R = (L_R, I_R, \Delta_\mathrm{re}(R))$ be a root system, namely, the tuple consisting of a free $\ZZ$-module $L_R$ of rank $\mu$, a symmetric $\ZZ$-bilinear form (called the {\em Cartan form}) $I_R \colon L_R \times L_R \longrightarrow \ZZ$ and a set of {\em real roots} $\Delta_\mathrm{re}(R) \subset L$ satisfying the following properties:
\begin{enumerate}
\item $L_R= \ZZ \Delta_\mathrm{re} (R) \cong \ZZ^\mu$,
\item $I_R(\alpha, \alpha) = 2$ for all $\alpha \in \Delta_\mathrm{re}(R)$,
\item Denote by $\Aut_\ZZ(L_R, I_R)$ the group of automorphisms of $L_R$ preserving the Cartan form $I_R$.
For every $\alpha \in \Delta_\mathrm{re}(R)$, the {\em reflection} $r_\alpha \in \Aut_\ZZ(L_R, I_R)$ of $\alpha$ defined by 
\begin{equation*}
r_\alpha(\lambda) = \lambda - I_R(\alpha, \lambda) \alpha, \quad \lambda \in L_R
\end{equation*}
satisfies $r_\alpha(\Delta_\mathrm{re}(R)) = \Delta_\mathrm{re}(R)$.
\end{enumerate}
Denote by the same symbol $I_R$ the natural extension of the Cartan form $I_R$ on $L_R \otimes_\ZZ \RR$.
For the Cartan form $I_R$ on $L_R \otimes_\ZZ \RR$, let $\mu_+, \mu_0, \mu_-$ denote the numbers of positive, $0$, negative eigenvalues of $I_R$, respectively.
We call $(\mu_+, \mu_0, \mu_-)$ the {\em signature} of the root system $R$.
A root system is called {\em finite} if its signature is $(\mu, 0 ,0)$. 
It is known that an irreducible root system is finite if and only if it is isomorphic to a root system of ADE type.

Denote by $\rad{I_R} \coloneqq \{ \lambda \in L_R \mid I_R(\lambda, \lambda') = 0,\ \lambda' \in L_R \}$ the radical of the root system $R$.
Then, the natural projection $L_R \longrightarrow L_R / \rad{I_R}$ induces a quotient root system $R_{/ \rad{I_R}} = (L_R', I_R', \Delta'_\mathrm{re}(R))$ such that $L_R' \coloneqq L_R / \rad{I_R}$, $I_R'$ is the symmetric $\ZZ$-bilinear form induced by $I_R$ and $\Delta'_\mathrm{re}(R)$ is the image of $\Delta_\mathrm{re}(R)$ by the projection.
A root system $R$ is called of {\em affine ADE type} if the signature of $R$ is $(\mu - 1, 1, 0)$ and the quotient root system $R_{/ \rad{I_R}}$ is isomorphic to a root system of ADE type.

For a root system $R = (L_R, I_R, \Delta_\mathrm{re}(R))$, the {\em Weyl group} $W_R$ is defined to be a subgroup of $\Aut_\ZZ(L_R, I_R)$ generated by reflections:
\begin{equation*}
W_R \coloneqq \br{r_\alpha \in \Aut_\ZZ(L_R, I_R) \mid \alpha \in \Delta_\mathrm{re}(R)}.
\end{equation*}
Put $\h_R \coloneqq \Hom_\ZZ (L_R, \CC)$ and $\h_R^\ast \coloneqq L_R \otimes_\ZZ \CC$.
Since there is a natural pairing
\begin{equation*}
\br{-,-} \colon \h_R^\ast \times \h_R \longrightarrow \CC, \quad (\lambda, x) \mapsto \br{\lambda, x} \coloneqq x(\lambda),
\end{equation*}
the Weyl group $W_R$ acts on $\h_R$ as 
\begin{equation*}
\br{\lambda, w(x)} \coloneqq \br{w^{-1}(\lambda), x}
\end{equation*}
for $\lambda \in \h_R^\ast, x \in \h_R$ and $w \in W_R$.
There is also a $\CC^\ast$-action on $\h_R$ defined by $(q \cdot x) (\lambda) \coloneqq q \cdot x(\lambda)$ for $\lambda \in \h_R^\ast, x \in \h_R$ and $q \in \CC^\ast$.
Note that the $\CC^\ast$-action on $\h_R \setminus \{0 \}$ is free.
\subsection{Generalized root system}
Let $R = (L_R, I_R, \Delta_\mathrm{re}(R))$ be a root system.
A subset $B = \{ \alpha_1, \dots, \alpha_\mu \}$ of $\Delta_\mathrm{re}(R)$ is called a {\em root basis} of $R$ if $\Delta_\mathrm{re}(R) = \{ w(\alpha_i) \in \Delta_\mathrm{re} (R) \mid w \in W(B), \alpha_i \in B \}$, where $W(B) = \br{r_{\alpha_1}, \dots, r_{\alpha_\mu}}$.
A root basis $B = \{ \alpha_1, \dots, \alpha_\mu \}$ is a $\ZZ$-basis of $L_R$ and satisfies $W(B) = W_R$.
The {\em Coxeter--Dynkin diagram} of a root basis $B = \{ \alpha_1, \dots, \alpha_\mu \}$ is defined as follows:
\begin{itemize}
\item The set of vertices is $B$.
\item The edge between the vertices $\alpha_i$ and $\alpha_j$ is given by 
\begin{eqnarray*}
\xymatrix{\circ_{\alpha_i} & \circ_{\alpha_j}} & \text{if} & I_R(\alpha_i, \alpha_j) = 0, \\
\xymatrix{\circ_{\alpha_i} \ar@{-}[r] & \circ_{\alpha_j}} & \text{if} & I_R(\alpha_i, \alpha_j) = -1, \\
\xymatrix{\circ_{\alpha_i} \ar@{-}[r]^{m} & \circ_{\alpha_j}} & \text{if} & I_R(\alpha_i, \alpha_j) = -m, \ (m \ge 2)\\
\xymatrix{\circ_{\alpha_i} \ar@{--}[r] & \circ_{\alpha_j}} & \text{if} & I_R(\alpha_i, \alpha_j) = 1, \\
\xymatrix{\circ_{\alpha_i} \ar@{==}[r] & \circ_{\alpha_j}} & \text{if} & I_R(\alpha_i, \alpha_j) = 2, \\
\xymatrix{\circ_{\alpha_i} \ar@{--}[r]^m & \circ_{\alpha_j}} & \text{if} & I_R(\alpha_i, \alpha_j) = m, \ (m \ge 3) . 
\end{eqnarray*}
\end{itemize}

An element $\mathbf{c} \in W_R$ is called a {\em Coxeter transformation} if there exists a root basis $B = \{ \alpha_1, \dots, \alpha_\mu \}$ such that $\mathbf{c} \coloneqq r_{\alpha_1} \cdots r_{\alpha_\mu} \in W_R$.
\begin{defn}[{\cite{S1,S3}}]
A {\em generalized root system} $(R, \mathbf{c})$ is a pair of a root system $R$ and a Coxeter transformation $\mathbf{c}$.
\end{defn}

Next, we consider the Euler form for a generalized root system.
\begin{defn}[{\cite[Definition 2.15]{NST}}]\label{defn : Euler form}
Let $(R, \mathbf{c})$ be a generalized root system.
An {\em Euler form} is a non-degenerated $\ZZ$-bilinear form $\chi \colon L_R \times L_R \longrightarrow \ZZ$ such that 
\begin{subequations}\label{eq : Euler}
\begin{enumerate}
\item We have $I_R = \chi + \chi^T$, namely, 
\begin{equation}\label{eq : Euler 1}
I_R(\lambda, \lambda') = \chi(\lambda, \lambda') + \chi(\lambda', \lambda), \quad \lambda, \lambda' \in L_R.
\end{equation}
\item (Serre duality) We have $\mathbf{c} = - \chi^{-1} \chi^T$, namely, 
\begin{equation}
\chi(\lambda, \lambda') = - \chi(\lambda', \mathbf{c}(\lambda)), \quad \lambda, \lambda' \in L_R.
\end{equation}
\end{enumerate}
\end{subequations}
\end{defn}
The following Lemma states the existence of an Euler form.
\begin{lem}[{\cite[Lemma 2.12]{NST}}]\label{lem : Euler by NST}
Let $(R, \mathbf{c})$ be a generalized root system and $B = \{ \alpha_1, \dots, \alpha_\mu \}$ a root basis with $\mathbf{c} = r_{\alpha_1} \cdots r_{\alpha_\mu}$.
The $\ZZ$-bilinear form $\chi \colon L_R \times L_R \longrightarrow \ZZ$ define by 
\begin{equation*}
\chi(\alpha_i, \alpha_j) \coloneqq
\begin{cases}
I_R(\alpha_i, \alpha_j), & i < j, \\
1, & i = j, \\
0, & i > j.
\end{cases}
\end{equation*}
is an Euler form.
\qed
\end{lem}
Since equations \eqref{eq : Euler} imply $I_R = \chi (\mathrm{id} - \mathbf{c})$, the Euler form is unique if $I_R$ is non-degenerate.

\section{Twist automorphism and extended affine Weyl group}\label{sec : Twist automorphism and extended affine Weyl group}
In this section, we introduce the twist automorphism and prove that the Dubrovin--Zhang extended affine Weyl group is isomorphic to our (modified) extended affine Weyl group.
\subsection{Generalized root system of affine ADE type}
From now on, we only consider generalized root systems of affine ADE type.
Let $A = (a_1, a_2, a_3)$ be a tuple of positive integers. 
Define a rational number $\chi_A \in \QQ$ and a positive integer $\mu_A \in \ZZ_{\ge 1}$ by 
\begin{equation*}
\chi_A \coloneqq \dfrac{1}{a_1} + \dfrac{1}{a_2} + \dfrac{1}{a_3} - 1, \quad \mu_A \coloneqq a_1 + a_2 + a_3 - 1,
\end{equation*}
and put $\ell_A \coloneqq \mathrm{lcm}(a_1, a_2, a_3)$.
One can easily check that $\chi_A > 0$ if and only if $A = (1, p, q)$, $(2, 2, r)$, $(2, 3, 3)$, $(2, 3, 4)$ or $(2, 3, 5)$, where $p, q, r \in \ZZ_{\ge 1}$.
Throughout this paper, we always assume that $\chi_A > 0$.
For a tuple $A = (a_1, a_2, a_3)$, we associate a Dynkin diagram $\Gamma_A$ (as in Figure \ref{fig : Dynkin diagram}) as follows:
\begin{itemize}
\item The set of vertices is $V_A \coloneqq \{ \mathbf{1} \} \cup \{ (i,j) \mid i = 1, 2, 3, \ j = 1, \dots, a_i - 1 \}$.
\item There is exactly one edge between $\mathbf{1}$ and $(i,1)$, and $(i,j)$ and $(i,j+1)$.
Otherwise, there is no edge.
\end{itemize}
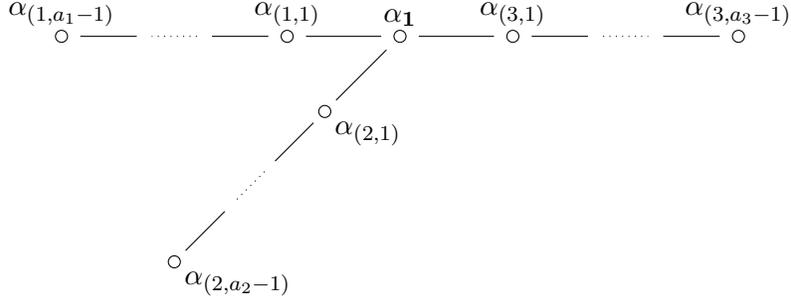
\begin{figure}[h]\label{fig : Coxeter-Dynkin of finite ADE}
\begin{tikzpicture}
\filldraw[fill=white] (0,0) circle[radius=0.8mm] node [above] {$\alpha_\mathbf{1}$};

\filldraw[fill=white] (-1.5,0) circle[radius=0.8mm] node [above] {$\alpha_{(1,1)}$};
\filldraw[fill=white] (-4.5,0) circle[radius=0.8mm] node [above] {$\alpha_{(1,a_1 - 1)}$};

\draw (-0.25,0) -- (-1.25,0);
\draw (-1.75,0) -- (-2.5,0);
\draw[dotted] (-2.7,0) -- (-3.3,0);
\draw (-3.5,0) -- (-4.25,0);

\filldraw[fill=white] (-1,-1) circle[radius=0.8mm] node [below right] {$\alpha_{(2,1)}$};
\filldraw[fill=white] (-3,-3) circle[radius=0.8mm] node [below right] {$\alpha_{(2,a_2 - 1)}$};

\draw (-0.18,-0.18) -- (-0.85,-0.85);
\draw (-1.15,-1.15) -- (-1.66,-1.66);
\draw[dotted] (-1.8,-1.8) -- (-2.2,-2.2);
\draw (-2.33,-2.33) -- (-2.85,-2.85);

\filldraw[fill=white] (1.5,0) circle[radius=0.8mm] node [above] {$\alpha_{(3,1)}$};
\filldraw[fill=white] (4.5,0) circle[radius=0.8mm] node [above] {$\alpha_{(3,a_3 - 1)}$};

\draw (0.25,0) -- (1.25,0);
\draw (1.75,0) -- (2.5,0);
\draw[dotted] (2.7,0) -- (3.3,0);
\draw (3.5,0) -- (4.25,0);
\end{tikzpicture}
\caption{Dynkin diagram $\Gamma_A$}
\label{fig : Dynkin diagram}
\end{figure}
The Dynkin diagram $\Gamma_A$ of $A = (1, p, q)$, $(2, 2, r)$, $(2, 3, 3)$, $(2, 3, 4)$ and $(2, 3, 5)$ are of type $A_{p + q - 1}$, $D_{r + 2}$, $E_6$, $E_7$ and $E_8$, respectively.

Define a generalized root system $(\widetilde{R}_A, \widetilde{\mathbf{c}}_A)$ of type $A$ as follows:
\begin{itemize}
\item $\widetilde{R}_A = (\widetilde{L}, \widetilde{I}, \widetilde{\Delta}_\mathrm{re})$ is a root system of affine ADE type whose rank is $\mu_A$, and the quotient root system $R_A \coloneqq (\widetilde{R}_{A})_{/ \rad{\widetilde{I}}} = (L, I, \Delta_\mathrm{re})$ admits a root basis $B_A = \{ \alpha_v \}_{v \in V_A}$ whose Coxeter--Dynkin diagram is the Dynkin diagram $\Gamma_A$.
\item Let $\delta \in \widetilde{L}$ be a generator of the radical $\rad{\widetilde{I}}$.
We regard $B_A$ as a subset of $\widetilde{L}$ by the isomorphism of lattices $(\widetilde{L}, \widetilde{I}) \cong (L, I) \oplus (\mathrm{rad}(\widetilde{I}), 0)$. 
Then, the set $\widetilde{B}_A = \{ \alpha_v \}_{v \in \widetilde{V}_A}$ is a root basis of $\widetilde{R}_A$, where $\widetilde{V}_A \coloneqq \{ \mathbf{1}^\ast \} \cup V_A$ and $\alpha_{\mathbf{1}^\ast} \coloneqq \alpha_\mathbf{1} + \delta$.
The coxeter transformation is defined by 
\begin{equation}\label{eq : affine Coxeter transformation}
\widetilde{\mathbf{c}}_A \coloneqq r_{\alpha_{\mathbf{1}^\ast}} r_{\alpha_\mathbf{1}} \prod_{i = 1}^3 \prod_{j = 1}^{a_i - 1} r_{\alpha_{(i, j)}}.
\end{equation}
\end{itemize}
For simplicity, we put $\widetilde{W}_A \coloneqq W_{\widetilde{R}_A}$, $W_A \coloneqq W_{R_A}$, $\widetilde{\h}_A \coloneqq \h_{\widetilde{R}_A}$, $\h_A \coloneqq \h_{R_A}$ and so on.
It is well-known that there is an isomorphism $\widetilde{W}_A \cong L \rtimes W_A$.

\begin{figure}[h]
\begin{tikzpicture}
\filldraw[fill=white] (0,0) circle[radius=0.8mm] node [below right] {$\alpha_\mathbf{1}$};
\filldraw[fill=white] (0,1.45) circle[radius=0.8mm] node [above right] {$\alpha_{\mathbf{1}^\ast}$};

\draw[dashed] (0.05,0.25) -- (0.05,1.25);
\draw[dashed] (-0.05,0.25) -- (-0.05,1.25);

\draw (-0.95,-0.82) -- (-0.1,1.3);
\draw (-1.35,0.15) -- (-0.15,1.35);
\draw (1.35,0.15) -- (0.15,1.35);

\filldraw[fill=white] (-1.5,0) circle[radius=0.8mm] node [below] {$\alpha_{(1,1)}$};
\filldraw[fill=white] (-4.5,0) circle[radius=0.8mm] node [below] {$\alpha_{(1,a_1 - 1)}$};

\draw (-0.25,0) -- (-1.25,0);
\draw (-1.75,0) -- (-2.5,0);
\draw[dotted] (-2.7,0) -- (-3.3,0);
\draw (-3.5,0) -- (-4.25,0);

\filldraw[fill=white] (-1,-1) circle[radius=0.8mm] node [below right] {$\alpha_{(2,1)}$};
\filldraw[fill=white] (-3,-3) circle[radius=0.8mm] node [below right] {$\alpha_{(2,a_2 - 1)}$};

\draw (-0.18,-0.18) -- (-0.85,-0.85);
\draw (-1.15,-1.15) -- (-1.66,-1.66);
\draw[dotted] (-1.8,-1.8) -- (-2.2,-2.2);
\draw (-2.33,-2.33) -- (-2.85,-2.85);

\filldraw[fill=white] (1.5,0) circle[radius=0.8mm] node [below] {$\alpha_{(3,1)}$};
\filldraw[fill=white] (4.5,0) circle[radius=0.8mm] node [below] {$\alpha_{(3,a_3 - 1)}$};

\draw (0.25,0) -- (1.25,0);
\draw (1.75,0) -- (2.5,0);
\draw[dotted] (2.7,0) -- (3.3,0);
\draw (3.5,0) -- (4.25,0);
\end{tikzpicture}
\caption{The Coxeter--Dynkin diagram $\widetilde{\Gamma}_A$}
\label{fig : Coxeter Dynkin diagram associated with the octopus quiver}
\end{figure}
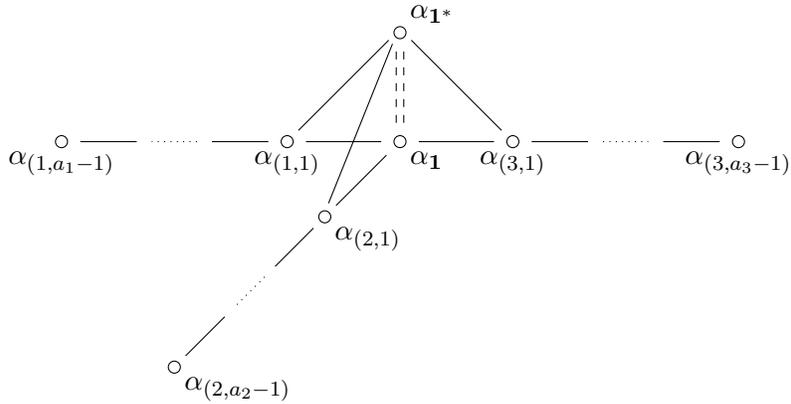

Note that the inclusion $B_A \subset \widetilde{B}_A$ implies a section $R_A \longrightarrow \widetilde{R}_A$.
In particular, we have 
\begin{equation}\label{eq : affine real roots}
\widetilde{\Delta}_\mathrm{re} = \{ \alpha + n \delta \in \widetilde{L} \mid \alpha \in \Delta_\mathrm{re}, n \in \ZZ \}.
\end{equation}
Denote by $\widetilde{\Gamma}_A$ the Coxeter--Dynkin diagram of $\widetilde{B}_A$. 
This diagram is given as in Figure \ref{fig : Coxeter Dynkin diagram associated with the octopus quiver}.

Denote by $\omega_v^\vee \in \h_A$ the fundamental co-weight of $\alpha_v \in B_A$, namely $\br{\alpha_{v'}, \omega_v^\vee} = \delta_{v' v}$, and $\alpha_v^\vee$ the co-root of $\alpha_v \in B_A$, namely $\br{\alpha_{v'}, \alpha_v^\vee} = I(\alpha_{v'}, \alpha_v)$.
Let $\delta^\vee \in \widetilde{\h}_A$ be the dual element to $\delta \in \widetilde{L}$.
Then,  
\begin{equation*}
\widetilde{\h}_A = \bigoplus_{v \in V_A} \CC \omega_v^\vee \oplus \CC \delta^\vee \cong \h_A \oplus \CC \delta^\vee.
\end{equation*}
It is known that the isomorphism 
\begin{equation}\label{eq : isomorphism and affine Weyl group}
\big\{ \widetilde{x} \in \widetilde{\h}_A \mid \br{\delta, \widetilde{x}} = 1 \big\} \overset{\cong}{\longrightarrow} \h_A, \quad
\widetilde{x} = x + \delta^\vee \mapsto x,
\end{equation}\label{eq : affine Weyl group equivariant isomorphism}
is $\widetilde{W}_A$-equivariant, where the group $\widetilde{W}_A \cong L \rtimes W_A$ acts on $\h_A$ by 
\begin{equation*}
x \mapsto w(x) + \sum_{v \in V_A} m_v \alpha_v^\vee, \quad w \in W_A, \ m_v \in \ZZ, \\
\end{equation*}
\subsection{Twist automorphism}
We define the twist automorphism for a generalized root system of affine ADE type.
At the first, we prove the uniqueness of the Euler form.
\begin{lem}[{cf.~\cite{NST}}]\label{lem : uniqueness of the Euler from}
Let $(\widetilde{R}_A, \widetilde{\mathbf{c}}_A)$ be a generalized root system of affine ADE type.
An Euler form $\chi \colon \widetilde{L} \times \widetilde{L} \longrightarrow \ZZ$ is unique.
\end{lem}
\begin{pf}
The existence follows from Lemma \ref{lem : Euler by NST}.
We prove the uniqueness here.
Let $\chi \colon \widetilde{L} \times \widetilde{L} \longrightarrow \ZZ$ be an Euler form of $(\widetilde{R}, \widetilde{\mathbf{c}}_A)$.
By the decomposition of root lattice $(\widetilde{L}, \widetilde{I}) \cong (L, I) \oplus (\mathrm{rad}(\widetilde{I}), 0)$, the restriction $\chi|_L \colon L \times L \longrightarrow \ZZ$ is uniquely determined.
It is sufficient to determine the values $\chi(\delta, \alpha_v)$ for $\alpha_v \in B_A$.
Indeed, by definition, we have $\chi(\delta, \delta) = 0$ and $\chi(\alpha_v, \delta) = - \chi(\delta, \alpha_v)$ for all $v \in V_A$.

By the definition \eqref{eq : affine Coxeter transformation} of the Coxeter transformation $\widetilde{\mathbf{c}}_A$, it follows form direct calculations that 
\begin{equation*}
\widetilde{\mathbf{c}}_A (\alpha_v) = 
\begin{cases}
\alpha_{\mathbf{1}} + \alpha_{(1,1)} + \alpha_{(2,1)} + \alpha_{(3,1)}, & v = \mathbf{1}^\ast, \\
\alpha_{\mathbf{1}} - \delta + \alpha_{(1,1)} + \alpha_{(2,1)} + \alpha_{(3,1)}, & v = \mathbf{1}, \\
\alpha_{(i, j + 1)}, & v = (i, j), \ i = 1, 2, 3, \ j = 1, \dots, a_i - 2, \\
\delta - \alpha_{(i, 1)} - \cdots - \alpha_{(i, a_i - 1)}, & v = (i, a_i - 1), \ i = 1, 2, 3.
\end{cases}
\end{equation*}
For $i = 1, 2, 3$ and $j = 1, \dots, a_i - 2$, we have $\chi(\delta, \alpha_{(i, j)}) = - \chi(\alpha_{(i, j)}, \delta) = \chi(\delta, \alpha_{(i, j + 1)})$.
If $j = a_i - 1$, we have 
\begin{eqnarray*}
\chi(\delta, \alpha_{(i, a_1 - 1)}) & = & \chi(\delta, \delta) - \sum_{j = 1}^{a_i - 1} \chi(\delta, \alpha_{(i,j)}) \\
& = & - \sum_{j = 1}^{a_i - 1} \chi(\delta, \alpha_{(i,j)}) \\
& = & - (a_i - 1) \chi(\delta, \alpha_{(i,a_i - 1)}),
\end{eqnarray*}
which yields $\chi(\delta, \alpha_v) = 0$ for $v \in V_A \setminus \{ \mathbf{1} \}$,

We calculate $\chi(\delta, \alpha_{\mathbf{1}})$.
Since $\chi|_L$ is unique, Lemma \ref{lem : Euler by NST} and Figure \ref{fig : Dynkin diagram} imply that 
\begin{eqnarray*}
\chi(\alpha_{\mathbf{1}^\ast}, \alpha_{\mathbf{1}}) & = & - \chi(\alpha_{\mathbf{1}}, \widetilde{\mathbf{c}}_A (\alpha_{\mathbf{1}^\ast})) \\
& = & - \chi|_L (\alpha_{\mathbf{1}}, \alpha_{\mathbf{1}}) - \chi|_L (\alpha_{\mathbf{1}}, \alpha_{(1,1)}) - \chi|_L (\alpha_{\mathbf{1}}, \alpha_{(2,1)}) - \chi|_L (\alpha_{\mathbf{1}}, \alpha_{(3,1)}) \\
& = & 2.
\end{eqnarray*}
Hence, $\chi(\delta, \alpha_{\mathbf{1}}) = \chi(\alpha_{\mathbf{1}^\ast}, \alpha_\mathbf{1}) - \chi (\alpha_{\mathbf{1}}, \alpha_\mathbf{1}) = 2 - 1 = 1$.
\qed
\end{pf}

\begin{defn}
Define an automorphism $t_\delta \in \Aut_\ZZ(\widetilde{L}, \widetilde{I})$ by 
\begin{equation}
t_\delta(\lambda) \coloneqq \lambda - \chi(\delta, \lambda) \delta, \quad \lambda \in \widetilde{L}.
\end{equation}
We call the automorphism the {\em twist automorphism} about $\delta$.
\end{defn}
One can easily check that $t_{\delta} = t_{- \delta}$.
Hence, the twist automorphism is independent of the choice of a generator of $\rad{\widetilde{I}}$.

\begin{rem}
From the view point of mirror symmetry, the twist automorphism is related to Dehn twists of a Riemannian surface.
One can construct a generalized root system of affine $A$ type from an annulus with marked points (cf.~\cite{HKK, IQ, STW}).
Under the construction, the Dehn twist about a unique simple closed curve induces the twist automorphism.
\end{rem}

We show the basic properties of the twist automorphism.
\begin{lem}\label{lem : basic properties}
Let $(\widetilde{R}_A, \widetilde{\mathbf{c}}_A)$ be a generalized root system of affine ADE type.
\begin{enumerate}
\item We have $t_\delta(\delta) = \delta$.
\item For $n \in \ZZ$, we have $t^n_\delta(\lambda) = \lambda - n \chi(\delta, \lambda) \delta$ for $\lambda \in \widetilde{L}$.
\item We have $t_\delta(\widetilde{\Delta}_\mathrm{re}) = \widetilde{\Delta}_\mathrm{re}$.
\item We have $t_\delta r_\alpha t_\delta^{-1} = r_{t_\delta(\alpha)}$ for $\alpha \in \widetilde{\Delta}_\mathrm{re}$.
In particular, for the root basis $\widetilde{B}_A = \{ \alpha_v \}_{v \in \widetilde{V}_A}$, we have 
\begin{equation*}
t_\delta r_{\alpha_v} t_\delta^{-1} = 
\begin{cases}
r_{\alpha_\mathbf{1}}, & v = \mathbf{1}^\ast, \\
r_{\alpha_{\mathbf{1}}} r_{\alpha_\mathbf{1}^\ast} r_{\alpha_{\mathbf{1}}}, & v = \mathbf{1}, \\
r_{\alpha_v}, & v \in \widetilde{V}_A \setminus \{ \mathbf{1}, \mathbf{1}^\ast \}.
\end{cases}
\end{equation*}
\end{enumerate}
\end{lem}
\begin{pf}
(1): By \eqref{eq : Euler 1}, we have $\chi(\delta, \delta) = 0$, which implies the statement.

(2): We prove the statement for $n \in \ZZ_{\ge 0}$.
By induction on $n$, it holds that 
\begin{eqnarray*}
t^n_\delta(\lambda) & = & t^{n - 1}_\delta(\lambda) - \chi(\delta, t^{n - 1}_\delta(\lambda)) \delta \\
& = & \lambda - (n - 1) \chi(\delta, \lambda) \delta - \chi(\delta, \lambda - (n - 1) \chi(\delta, \lambda) \delta) \delta \\
& = & \lambda - (n - 1) \chi(\delta, \lambda) \delta - \chi(\delta, \lambda) \delta \\
& = & \lambda - n \chi(\delta, \lambda) \delta,
\end{eqnarray*}
for any $\lambda \in \widetilde{L}$.
Similarly, one can show the statement for $n \in \ZZ_{\le 0}$.

(3): The statement follows from \eqref{eq : affine real roots}.

(4): For any $\lambda \in \widetilde{L}$, it holds that 
\begin{eqnarray*}
t_\delta r_\alpha t_\delta^{-1} (\lambda) 
& = & \lambda + \chi(\delta, \lambda) \delta - I(\alpha, \lambda) \alpha - \chi(\delta, \lambda + \chi(\delta, \lambda) \delta - I(\alpha, \lambda) \alpha) \delta \\
& = & \lambda - I(\alpha - \chi(\delta, \alpha) \delta, \lambda) (\alpha - \chi(\delta, \alpha) \delta) \\
& = & r_{t_\delta(\alpha)} (\lambda).
\end{eqnarray*}

By the proof of Lemma \ref{lem : uniqueness of the Euler from}, we have 
\begin{equation*}
t_\delta (\alpha_v) = 
\begin{cases}
\alpha_\mathbf{1}, & v = \mathbf{1}^\ast, \\
\alpha_\mathbf{1} - \delta, & v = \mathbf{1}, \\
\alpha_v, & v \in \widetilde{V}_A \setminus \{ \mathbf{1}, \mathbf{1}^\ast \}.
\end{cases}
\end{equation*}
Hence, a direct calculation implies the last statement.
\qed
\end{pf}
\subsection{Extended affine Weyl groups}
In this section, we define an extended affine Weyl group and a modified one.

\begin{defn}\label{defn : extended affine Weyl group}
Let $(\widetilde{R}_A, \widetilde{\mathbf{c}}_A)$ be a generalized root system of affine ADE type.
\begin{enumerate}
\item Define a group $T_A \coloneqq \br{t_\delta}$ as the subgroup of $\Aut_\ZZ (\widetilde{L}, \widetilde{I})$ generated by the twist automorphism.
\item Define a group $\widetilde{W}_A^\mathrm{ext}$ as the subgroup of $\Aut_\ZZ (\widetilde{L}, \widetilde{I})$ generated by the affine Weyl group $\widetilde{W}_A$ and $T_A$.
We call $\widetilde{W}_A^\mathrm{ext}$ the {\em extended affine Weyl group}.
\end{enumerate}
\end{defn}
By Lemma \ref{lem : basic properties} (2), there is an isomorphism $T_A \cong \ZZ$.

\begin{cor}
There exists a short exact sequence of groups
\begin{equation}\label{eq : ses for W and T}
1 \longrightarrow \widetilde{W}_A \longrightarrow \widetilde{W}_A^\mathrm{ext} \overset{p}{\longrightarrow} T_A \longrightarrow 1.
\end{equation}
\end{cor}
\begin{pf}
The statement follows from Lemma \ref{lem : basic properties} (3) and (4).
\qed
\end{pf}

For an element $\lambda \in \widetilde{\h}_A^\ast$, denote by $H_\lambda \coloneqq \{ \widetilde{x} \in \widetilde{\h}_A \mid \br{\lambda, \widetilde{x}} = 0 \}$ the complex orthogonal hyperplane. 
We consider subspaces $X_A$ and $X_A^\mathrm{reg}$, which was considered by Bridgeland \cite[Section 1.3]{Br2} and Ikeda \cite[Definition 2.7]{I} in the study of the space of Bridgeland stability conditions:
\begin{enumerate}
\item Define a subset $X_A \subset \widetilde{\h}_A$ by
\begin{equation*}
X_A \coloneqq \widetilde{\h}_A \setminus H_\delta.
\end{equation*}
\item Define the {\em regular subspace} $X_A^\mathrm{reg} \subset X_A$ by 
\begin{equation*}
X_A^\mathrm{reg} \coloneqq X_A \setminus \bigcup_{\alpha \in \widetilde{\Delta}_{\rm re}} H_\alpha.
\end{equation*}
\end{enumerate}
In order to study Frobenius manifold associated with a generalized root system, we consider the following spaces:
\begin{enumerate}
\item[(3)] The universal covering space $\pi_X \colon \X_A \longrightarrow X_A$ of $X_A$,
\item[(4)] The {\em regular subspace} $\X_A^\mathrm{reg} \coloneqq \pi_X^{-1}(X_A^\mathrm{reg})$.
\end{enumerate}

We define a $T_A$-action on $\widetilde{\h}_A$ by $\br{\lambda, t_\delta(\widetilde{x})} = \br{t_\delta^{-1}(\lambda), \widetilde{x}}$ for $\lambda \in \widetilde{\h}_A^\ast$ and $\widetilde{x} \in \widetilde{\h}_A$.
This action also defines the $\widetilde{W}_A^\mathrm{ext}$-action on $\widetilde{\h}_A$.
By Lemma \ref{lem : basic properties} (1), the $\widetilde{W}_A^\mathrm{ext}$-action on $\widetilde{\h}_A$ induces the one on $X_A$.
Moreover, Lemma \ref{lem : basic properties} (3) implies that the action induces the one on $X_A^\mathrm{reg}$.
We also have the free $\CC^\ast$-action on $X_A$ and $X_A^\mathrm{reg}$.
The $\widetilde{W}_A^\mathrm{ext}$-action and the $\CC^\ast$-action commute.

\begin{lem}\label{lem : isomorphism between X and h times C}
The map 
\begin{equation*}
\varphi \colon X_A \longrightarrow \h_A \times \CC^\ast, \quad 
\sum_{v \in V_A} x_v \omega_v^\vee + x_\delta \delta^\vee \mapsto \Big( \dfrac{1}{x_\delta} \sum_{v \in V_A} x_v \omega_v^\vee, x_\delta \Big), 
\end{equation*}
is a $\widetilde{W}_A^\mathrm{ext}$-equivariant biholomorphism, where the $\widetilde{W}_A^\mathrm{ext}$-action on $\h_A \times \CC^\ast$ is given as follows: \\
The action of $\widetilde{W}_A \cong L \rtimes W_A$ is 
\begin{equation*}
(x, z) \mapsto \Big( w(x) + \sum_{v \in V_A} m_v \alpha_v^\vee, z \Big), \quad w \in W_A, \ m_v \in \ZZ.
\end{equation*}
The action of $T_A \cong \ZZ$ is 
\begin{equation*}
(x, z) \mapsto (x + m \omega_{\mathbf{1}}^\vee, z), \quad m \in \ZZ \cong T.
\end{equation*}
In particular, we have $\pi_1(X_A) \cong \ZZ$.
\end{lem}
\begin{pf}
By definition, it is clear that the map $\varphi$ is a biholomorphism.
Since the $\CC^\ast$-action is free, it follows from \eqref{eq : isomorphism and affine Weyl group} that the map $\varphi$ is $\widetilde{W}_A$-equivariant.

Let $\lambda  = \sum_{v \in V_A} \lambda_v \alpha_v + \lambda_\delta \delta \in \widetilde{\h}_A^\ast$.
By Lemma \ref{lem : uniqueness of the Euler from}, it holds that
\begin{equation*}
\chi(\delta, \lambda) = \sum_{v \in V_A} \lambda_v \chi(\delta, \alpha_v) = \lambda_\mathbf{1} = \br{\lambda, \omega_\mathbf{1}^\vee}.
\end{equation*}
Therefore, for $\widetilde{x} = \sum_{v \in V_A} x_v \omega_v^\vee + x_\delta \delta^\vee \in X_A$ and $\lambda \in \widetilde{\h}_A^\ast$, we have
\begin{eqnarray*}
\br{\lambda, t_\delta(\widetilde{x})} & = & \br{\lambda + \chi(\delta, \lambda) \delta, \widetilde{x}} \\
& = & \br{\lambda, \widetilde{x}} + \br{\lambda, \omega_\mathbf{1}^\vee} \br{\delta, \widetilde{x}} \\
& = & \br{\lambda, \widetilde{x} + \br{\delta, \widetilde{x}} \omega_\mathbf{1}^\vee}.
\end{eqnarray*}
Thus, we obtain 
\begin{equation*}
\varphi (t_\delta (\widetilde{x})) = \varphi( \widetilde{x} + \br{\delta, \widetilde{x}} \omega_\mathbf{1}^\vee) 
= \Big( \dfrac{x_\mathbf{1} + x_\delta}{x_\delta} \omega_\mathbf{1}^\vee + \sum_{v \ne \mathbf{1}} \dfrac{x_v}{x_\delta} \omega_v^\vee, x_\delta \Big) 
= t_\delta(\varphi(\widetilde{x})).
\end{equation*}
\qed
\end{pf}
By Lemma \ref{lem : isomorphism between X and h times C}, one can define a $\widetilde{W}_A^\mathrm{ext}$-action on $\X_A$.
However, in order to construct Frobenius manifold for the generalized affine root system $(\widetilde{R}_A, \widetilde{\mathbf{c}}_A)$, it is necessary to modify the group $\widetilde{W}_A^\mathrm{ext}$ by the fundamental group $\pi_1(X_A)$.
Let $\gamma$ be the generator of the fundamental group $\pi_1(X_A) \cong \ZZ$ corresponding to $1 \in \ZZ$.
Denote by $C_A = \br{(\gamma, t_\delta)}$ the diagonal subgroup of $T_A \times \pi_1(X_A)$ generated by the element $(\gamma, t_\delta)$.
Then there is a canonical isomorphism $\phi \colon C_A \cong T_A$.
\begin{equation*}
\xymatrix{
& T_A \times \pi_1(X_A) \ar[ld] \ar[rd] & \\ 
T_A & C_A \ar[l]^-{\phi}_-\cong \ar[r]^-\cong \ar@{}[u]|{\bigcup} & \pi_1(X_A) 
}
\end{equation*}
\begin{defn}
Let $p \colon \widetilde{W}_A^\mathrm{ext} \longrightarrow T_A$ be the natural projection in \eqref{eq : ses for W and T}.
Define a group $\widehat{W}_A$ by the fiber product of $C_A$ with respect to the map $p \colon \widetilde{W}_A^\mathrm{ext} \longrightarrow T_A$ and the map $\phi \colon C_A \cong T_A$:
\begin{equation}
\widehat{W}_A \coloneqq \widetilde{W}_A^\mathrm{ext} \times_{T_A} C_A = \{ (w, g) \in \widetilde{W}_A^\mathrm{ext} \times C_A \mid p(w) = \phi(g) \}
\end{equation}
We call the group $\widehat{W}_A$ the {\em modified extended affine Weyl group}.
\end{defn}
By the universal property of the fiber product, there exists a following commutative diagram whose rows are exact sequences:
\begin{equation*}
\xymatrix{
1 \ar[r] & \widetilde{W}_A \ar[r] \ar@{=}[d] & \widehat{W}_A \ar[r] \ar[d] & C_A \ar[r] \ar[d]^{\cong}_{\phi} & 1 \\
1 \ar[r] & \widetilde{W}_A \ar[r] & \widetilde{W}_A^\mathrm{ext} \ar[r]^{p} & T_A \ar[r] & 1
}
\end{equation*}
The five lemma yields a group isomorphism $\widehat{W}_A \cong \widetilde{W}_A^\mathrm{ext}$.
However, note that the $\widehat{W}_A$-action and $\widetilde{W}_A^\mathrm{ext}$-action on $\X_A$ are different.

Here, we recall the definition of Dubrovin--Zhang extended affine Weyl group.
\begin{defn}[{\cite[p.171]{DZ}}]\label{defn : Dubrovin--Zhang extended affine Weyl group}
Let $\widehat{\h}_A \coloneqq \h_A \times \CC$ and $x_{\mu_A}$ the coordinate on $\CC$.
\begin{subequations}\label{eq : DZ}
Define an action $\widetilde{W}_A$ on $\widehat{\h}_A$ by 
\begin{equation}\label{eq : DZ1}
(x, x_{\mu_A}) \mapsto \Big( w(x) + \sum_{v \in V_A} m_v \alpha_v^\vee, x_{\mu_A} \Big), \quad w \in W_A, \ m_v \in \ZZ, 
\end{equation}
and an action $\ZZ$ on $\widehat{\h}_A$ by  
\begin{equation}\label{eq : DZ2}
(x, x_{\mu_A}) \mapsto (x + m \omega_\mathbf{1}^\vee, x_{\mu_A} + m), \quad m \in \ZZ.
\end{equation}
\end{subequations}
Define the {\em Dubrovin--Zhang extended affine Weyl group} $\widehat{W}_A^\mathrm{DZ}$ as the subgroup of $\Aut(\widehat{\h}_A)$ generated by the transformations \eqref{eq : DZ}.
\end{defn}

\begin{prop}\label{prop : extended affine Weyl groups}
There exists an isomorphism $\widehat{W}_A^\mathrm{DZ} \cong \widehat{W}_A$.
\end{prop}
\begin{pf}
The map $\widehat{\pi} \colon \widehat{\h}_A \longrightarrow \h_A \times \CC^\ast$ defined by $\widehat{\pi}(x, x_{\mu_A}) \coloneqq (x, e^{2 \pi \sqrt{-1} x_{\mu_A}})$ is the universal covering map, Lemma \ref{lem : isomorphism between X and h times C} implies that the group $\pi_1(X_A)$ acts on $\widehat{\h}_A$ by $\gamma \cdot (x, x_{\mu_A}) = (x, x_{\mu_A} + 1)$.
Since the $\widetilde{W}_A^\mathrm{ext}$-action on $\CC^\ast$ is trivial, there is a group action of $\widetilde{W}_A^\mathrm{ext} \times \pi_1(X_A)$ on $\widehat{\h}_A$.
On the other hand, the inclusion $C_A \subset T_A \times \pi_1(X_A)$ yields an inclusion $\widehat{W}_A \subset \widetilde{W}_A^\mathrm{ext} \times_{T_A} (T_A \times \pi_1(X_A)) \cong \widetilde{W}_A^\mathrm{ext} \times \pi_1(X_A)$.
Hence, there is a group action of $\widehat{W}_A$ on $\widehat{\h}_A$.
It follows from Lemma \ref{lem : isomorphism between X and h times C} that the $\widehat{W}_A$-action on $\widehat{\h}_A$ coincides with \eqref{eq : DZ}.
Since the group representation $\widehat{W}_A \longrightarrow \Aut(\widehat{\h}_A)$ is an injection, we obtain the statement.
\qed
\end{pf}
As a direct consequence, we obtain the following
\begin{cor}
There exists an isomorphism of groups $\widehat{W}_A^\mathrm{DZ} \cong \widetilde{W}_A^\mathrm{ext}$.
\qed
\end{cor}

The following theorem is important to obtain the Frobenius manifold constructed by Dubrovin--Zhang in terms of a generalized affine root system.
\begin{thm}\label{thm : main 1}
There exists a $\widehat{W}_A$-equivariant biholomorphism $\widehat{\varphi} \colon \X_A \longrightarrow \widehat{\h}_A$.
\end{thm}
\begin{pf}
By the universal property, there exists a biholomorphism $\widehat{\varphi} \colon \X_A \longrightarrow \widehat{\h}_A$ such that $\widehat{\pi} \circ \widehat{\varphi} = \varphi \circ \pi_X$.
Let $L(\widetilde{W}_A^\mathrm{ext}) \coloneqq \{ \widetilde{W}_A \in \Aut(\X_A) \mid \pi_X \circ \widetilde{W}_A = w \circ \pi_X \ \text{for some $w \in \widetilde{W}_A^\mathrm{ext}$} \}$ be the lift of the $\widetilde{W}_A^\mathrm{ext}$-action to $\X_A$.
For $\widetilde{w} \in L(\widetilde{W}_A^\mathrm{ext})$ and $w \in \widetilde{W}_A^\mathrm{ext}$, by Lemma \ref{lem : isomorphism between X and h times C} we have
\begin{eqnarray*}
\pi_X \circ \big( \widehat{\varphi}^{-1} \circ \widetilde{w} \circ \widehat{\varphi} \big) 
& = & \varphi^{-1} \circ \widehat{\pi} \circ \widetilde{w} \circ \widehat{\varphi} \\
& = & \varphi^{-1} \circ w \circ \widehat{\pi} \circ \widehat{\varphi} \\
& = & w \circ \varphi^{-1} \circ \widehat{\pi} \circ \widehat{\varphi} \\
& = & w \circ \pi_X.
\end{eqnarray*}
Hence the map 
\begin{equation*}
\pi_1(X_A) \times \widetilde{W}_A^\mathrm{ext} \longrightarrow L(\widetilde{W}_A^\mathrm{ext}), \quad \widetilde{w} \mapsto \widehat{\varphi}^{-1} \circ \widetilde{w} \circ \widehat{\varphi}
\end{equation*}
is well-defined and obviously an isomorphism.
Hence, the biholomorphism $\widehat{\varphi}$ is equivariant with respect to $\pi_1(X_A) \times \widetilde{W}_A^\mathrm{ext}$.
In particular, since $\widehat{W}_A \subset \widetilde{W}_A^\mathrm{ext} \times \pi_1(X_A)$, we obtain the statement.
\qed
\end{pf}

Let $\CC[\X_A]$ be the polynomial $\CC$-algebra over $\X_A$ and $\CC[\X_A]^{\widehat{W}_A}$ the $\widehat{W}_A$-invariant subalgebra.
Denote by $\X_A / \widehat{W}_A$ the complex manifold corresponding to the algebraic variety $\mathrm{Spec} \ \CC[\X_A]^{\widehat{W}_A}$.  
\begin{cor}\label{cor : orbit space}
There exists a biholomorphism $\X_A / \widehat{W}_A \cong \CC^{\mu_A - 1} \times \CC^\ast$.
\end{cor}
\begin{pf}
It follows from \cite[Theorem 1.1]{DZ} that $\widehat{\h}_A / \widehat{W}_A^\mathrm{DZ} \cong \CC^{\mu_A - 1} \times \CC^\ast$.
Theorem \ref{thm : main 1} implies the isomorphism $\X_A / \widehat{W}_A \cong \widehat{\h}_A / \widehat{W}_A^\mathrm{DZ}$.
\qed
\end{pf}

\begin{cor}\label{cor : Frobenius structure}
Let $M = \X_A / \widehat{W}_A$, $\O_M$ the structure sheaf and $\T_M$ the tangent sheaf.
The complex manifold $M$ has a Frobenius structure $(\eta, \circ, e, E)$ whose conformal dimension is $1$, where 
\begin{itemize}
\item $\eta \colon \T_M \times \T_M \longrightarrow \O_M$ is a symmetric non-degenerate $\O_M$-bilinear form,
\item $\circ \colon \T_M \times \T_M \longrightarrow \T_M$ is a commutative $\O_M$-bilinear product,
\item $e \in \Gamma(M, \T_M)$ is the unit vector field of the product, 
\item $E \in \Gamma(M, \T_M)$ is a vector field, called the Euler vector field. 
\end{itemize}
Namely, the tuple satisfies the following conditions:
\begin{enumerate}
\item The product $\circ$ is self-adjoint with respect to $\eta$.
\item The Levi--Civita connection $\ns \colon \T_M \otimes_{\O_M} \T_M \longrightarrow \T_M$ with respect to $\eta$ is flat.
\item The $\O_M$-linear morphism $C_{\xi} \colon \T_M \longrightarrow \T_M$ for $\xi \in \T_M$ defined by $C_{\xi}(-) \coloneqq \xi \circ -$ is $\ns$-flat.
\item The unit vector field $e$ is $\ns$-flat.
\item We have $\mathrm{Lie}_E (\eta) = \eta$ and $\mathrm{Lie}_E (\circ) = \circ$.
\end{enumerate}
\end{cor}
\begin{pf}
By Corollary \ref{cor : orbit space}, the statement follows directly from \cite[Theorem 2.1]{DZ}.
\qed
\end{pf}

We recall the Lyashko--Looijenga map of the Frobenius manifold $\X_A / \widehat{W}_A$ defined in \cite[Section 3]{DZ}.
It is known by \cite[Section 3]{DZ} that there exists a local coordinate system $(u_1, \dots, u_{\mu_A})$, called the {\em canonical coordinate system}, such that 
\begin{equation*}
e = \frac{\p}{\p u_1} + \dots + \frac{\p}{\p u_{\mu_A}}, \quad 
E = u_1 \frac{\p}{\p u_1} + \dots + u_{\mu_A} \frac{\p}{\p u_{\mu_A}}, 
\end{equation*}
\begin{equation*}
\frac{\p}{\p u_i} \circ \frac{\p}{\p u_j} = \delta_{ij} \frac{\p}{\p u_i}, \quad i, j = 1, \dots, \mu_A, 
\end{equation*}
where $\delta_{ij}$ is the Kronecker's delta.
The canonical coordinate system $(u_1, \dots, u_{\mu_A})$ is uniquely determined up to a permutation of indices.
Then, the {\em Lyashko--Looijenga map} is defined by
\begin{equation*}
LL \colon \X_A / \widehat{W}_A \longrightarrow \CC^{\mu_A}, \quad LL (t) \coloneqq (b_1, \dots, b_{\mu_A}),
\end{equation*}
where $b_1, \dots, b_{\mu_A}$ are coefficients of the characteristic polynomial of $C_E$:
\begin{equation*}
\det (C_{E} - w) = \prod_{i = 1}^{\mu_A} (u_i - w) = (-1)^{\mu_A} (w^{\mu_A}  + b_1 w^{\mu_A- 1} + \dots + b_{\mu_A}).
\end{equation*}
The Lyashko--Looijenga map $LL$ is a ramified covering map and its {\em degree} $\deg LL$ is given explicitly as follows:
\begin{equation*}
\deg LL = \frac{\mu_A!}{a_1!a_2!a_3!\chi_A}a_1^{a_1}a_2^{a_2}a_3^{a_3} .
\end{equation*}

Next, we recall the monodromy group of the Frobenius manifold.
The subset $\{ t \in \X_A / \widehat{W}_A \mid \det (C_E)_t = 0 \}$ is called {\em discriminant locus} of the Frobenius manifold.
By the construction of the Frobenius structure \cite[Theorem 2.1]{DZ}, Theorem \ref{thm : main 1} implies that the complement of the discriminant locus is given by $\X_A^\mathrm{reg} / \widehat{W}_A$.
The {\em intersection form} $g \colon \T_{\X_A^\mathrm{reg} / \widehat{W}_A} \times \T_{\X_A^\mathrm{reg} / \widehat{W}_A} \longrightarrow \O_{\X_A^\mathrm{reg} / \widehat{W}_A}$ of the Frobenius manifold $\X_A / \widehat{W}_A$ is defined by $g (\xi, \xi') \coloneqq \eta(C_E^{-1} \xi, \xi')$ for $\xi, \xi' \in \T_{\X_A^\mathrm{reg} / \widehat{W}_A}$.
Then, flat functions with respect to the Levi--Civita connection of the intersection form $g$ define a local system on $\X_A^\mathrm{reg} / \widehat{W}_A$.
The fundamental group $\pi_1(\X_A^\mathrm{reg} / \widehat{W}_A)$ acts on a fiber of the local system.
Since the fiber of the local system can be identified with $\X_A$, we obtain a representation
\begin{equation*}
\rho \colon \pi_1(\X_A^\mathrm{reg} / \widehat{W}_A) \longrightarrow \Aut(\X_A).
\end{equation*}
\begin{defn}[{\cite[Definition G.1]{D}}]\label{defn : monodromy of Frobenius manifold}
The {\em monodromy group} of the Frobenius manifold $\X_A / \widehat{W}_A$ is defined to be the image of the map $\rho$.
\end{defn}

\begin{cor}
The monodromy group of the Frobenius manifold $\X_A / \widehat{W}_A$ is isomorphic to the group $\widehat{W}_A$.
\end{cor}
\begin{pf}
The statement follows from \cite[Corollary 2.7]{DZ} and Corollary \ref{cor : orbit space}.
\qed
\end{pf}

\section{Twist automorphism and spherical twist functor}\label{sec : Twist automorphism and spherical twist functor}
In this section, we study the relation between the set of root bases with the Coxeter transformation $\widetilde{\mathbf{c}}_A$ and the set of full exceptional collections in the derived category of a certain algebra.

We first review spherical twists and full exceptional collections.
Let $\D$ be a $\CC$-linear triangulated category of finite type with the Serre functor $\S_\D \in \Aut(\D)$.
We assume that the triangulated category $\D$ is equivalent to the perfect derived category of a smooth proper differential graded $\CC$-algebra.
By this assumption, we can have $\RR \Hom$-complexes and functorial cones.

Let $\XX$ be a distinguished autoequivalence on $\D$.
It will be convenient to put $[n + m \XX] \coloneqq [n] \circ \XX^m$ for $n, m \in \ZZ$ and $\Hom^\bullet_{\D} (X, Y) \coloneqq \bigoplus_{n, m \in \ZZ} \Hom_{\D} (X[n + m \XX], Y) [n + m \XX]$.
\begin{defn}[{\cite[Section A.4]{IQ}, \cite[Definition 1.1]{ST}}]
An object $S \in \D$ is called {\em $\XX$-spherical} if $\S_\D(S)\cong S[\XX]$ and $\Hom^\bullet_{\D} (S, S) \cong \CC \oplus \CC[-\XX]$.
\end{defn}

\begin{prop}[{\cite[Proposition A.8]{IQ}, \cite[Proposition 2.10]{ST}}]
For an $\XX$-spherical object $S \in \D$, there exists an autoequivalence $\T_S \in \Aut(\D)$ defined by the exact triangle 
\begin{subequations}
\begin{equation}\label{eq : spherical twist}
\Hom^\bullet_{\D} (S, X) \otimes S \longrightarrow X \longrightarrow \T_S(X), \quad X \in \D.
\end{equation}
The inverse functor $\T_S^{-1} \in \Aut(\D)$ is given by 
\begin{equation}
\T_S^{-1}(X) \longrightarrow X \longrightarrow S \otimes \Hom^\bullet_{\D} (X, S)^\ast, \quad X \in \D,
\end{equation}
where $(-)^\ast$ denotes the duality $\Hom_\CC(-, \CC)$.
\end{subequations}
\qed
\end{prop}
\begin{prop}[{\cite[Proposition 2.6]{ST}}]\label{prop : spherical object and twist}
Let $S$ and $S'$ be $\XX$-spherical objects in $\D$.
Assume that $S \cong S'[n + m \XX]$ for some $n, m \in \ZZ$.
Then, we have $\T_S \cong \T_{S'}$.
\qed
\end{prop}

Next, we recall full exceptional collections.
\begin{itemize}
\item An object $E \in \D$ is called {\em exceptional} if $\Hom^\bullet_\D(E,E) \cong \CC$.
\item An ordered set $\E = (E_1, \dots, E_\mu)$ consisting of exceptional objects $E_1, \dots, E_\mu$ is called {\em exceptional collection} if $\Hom^\bullet_\D(E_i, E_j) \cong 0$ for $i > j$.
\item An exceptional collection $\E$ is called {\em full} if the smallest full triangulated subcategory of $\D$ containing all elements in $\E$ is equivalent to $\D$.
\item Two full exceptional collections $\E=(E_1, \dots, E_\mu)$, $\F=(F_1, \dots, F_\mu)$ in $\D$ are {\em isomorphic} if $E_i \cong F_i$ for all $i = 1, \dots, \mu$.
\end{itemize}
Denote by $\FEC(\D)$ the set of isomorphism classes of full exceptional collections in $\D$.
\subsection{Derived category and generalized affine root system}
In this section, we explain the relation between the twist automorphism and spherical twists on the derived category which defines a generalized affine root system.

\begin{defn}[{\cite[Definition 2.19]{STW}}]\label{defn : octopus algebra}
Let $A = (a_1, a_2, a_3)$ with $\chi_A > 0$.
Define a quiver $Q_{\widetilde{\TT}_A}$ as in Figure \ref{fig : octopus}.
More precisely, 
\begin{itemize}
\item The set of vertices is given by $\widetilde{V}_A$.
\item Let $v, v' \in \widetilde{V}_A$ be vertices.
\begin{itemize}
\setlength{\leftskip}{-10pt}
\item If $(v, v') = (\mathbf{1}, (a_i, 1))$ or $((a_i, 1), \mathbf{1}^*)$, there is one arrow $\alpha_{v, v'} \colon v \to v'$.
\item If $(v, v') = ((a_i, j), (a_i, j + 1))$, there is one arrow $\alpha_{v, v'} \colon v \to v'$.
\item Otherwise, there are no arrows.
\end{itemize}
\end{itemize}
Define an admissible ideal $\I_{\widetilde{\TT}_A}$ of the path algebra $\CC Q_{\widetilde{\TT}_A}$ by 
\begin{equation*}
\I_{\widetilde{\TT}_A} \coloneqq \langle \alpha_{(1, 1), \mathbf{1}^*} \alpha_{\mathbf{1}, (1, 1)} + \alpha_{(2, 1), \mathbf{1}^*} \alpha_{\mathbf{1}, (2, 1)}, \ \alpha_{(2, 1), \mathbf{1}^*} \alpha_{\mathbf{1}, (2, 1)} + \alpha_{(3, 1), \mathbf{1}^*} \alpha_{\mathbf{1}, (3, 1)} \rangle.
\end{equation*}
The algebra $\CC \widetilde{\TT}_A \coloneqq \CC Q_{\widetilde{\TT}_A} / \I_{\widetilde{\TT}_A}$ is called the {\em octopus algebra} of type $A$.
\end{defn}
\begin{figure}[h]
\begin{tikzpicture}
\filldraw[fill=white] (0,0) circle[radius=0.8mm] node [below right] {\footnotesize $\mathbf{1}$};
\filldraw[fill=white] (0,1.45) circle[radius=0.8mm] node [above right] {\footnotesize ${\mathbf{1}^\ast}$};

\draw [->] (-0.95,-0.82) -- (-0.05,1.3);
\draw [->] (-1.35,0.15) -- (-0.15,1.35);
\draw [->] (1.35,0.15) -- (0.15,1.35);

\filldraw[fill=white] (-1.5,0) circle[radius=0.8mm] node [below] {\footnotesize $(1,1)$};
\filldraw[fill=white] (-4.5,0) circle[radius=0.8mm] node [below] {\footnotesize $(1,a_1 - 1)$};

\draw [->] (-0.25,0) -- (-1.25,0);
\draw [->] (-1.75,0) -- (-2.5,0);
\draw [dotted] (-2.7,0) -- (-3.3,0);
\draw [->] (-3.5,0) -- (-4.25,0);

\filldraw[fill=white] (-1,-1) circle[radius=0.8mm] node [below right] {\footnotesize $(2,1)$};
\filldraw[fill=white] (-3,-3) circle[radius=0.8mm] node [below right] {\footnotesize $(2,a_2 - 1)$};

\draw [->] (-0.18,-0.18) -- (-0.85,-0.85);
\draw [->] (-1.15,-1.15) -- (-1.66,-1.66);
\draw [dotted] (-1.8,-1.8) -- (-2.2,-2.2);
\draw [->] (-2.33,-2.33) -- (-2.85,-2.85);

\filldraw[fill=white] (1.5,0) circle[radius=0.8mm] node [below] {\footnotesize $(3,1)$};
\filldraw[fill=white] (4.5,0) circle[radius=0.8mm] node [below] {\footnotesize $(3,a_3 - 1)$};

\draw [->] (0.25,0) -- (1.25,0);
\draw [->] (1.75,0) -- (2.5,0);
\draw [dotted] (2.7,0) -- (3.3,0);
\draw [->] (3.5,0) -- (4.25,0);
\end{tikzpicture}
\caption{Quiver $Q_{\widetilde{\TT}_A}$}
\label{fig : octopus}
\end{figure}
For simplicity, we denote by $\D_A \coloneqq \D^b\mod (\CC \widetilde{\TT}_A)$ the bounded derived category of finitely generated modules over the algebra $\CC \widetilde{\TT}_A$.

On the other hand, one can also associate an extended Dynkin quiver with $A$.
For $A = (1, p, q)$, $(2, 2, r)$, $(2, 3, 3)$, $(2, 3, 4)$ and $(2, 3, 5)$, let $Q_A$ denote the extended Dynkin quiver $A_{p, q}^{(1)}$, $D_{r + 2}^{(1)}$, $E_6^{(1)}$, $E_7^{(1)}$ and $E_8^{(1)}$, respectively.
Let $\D^b(Q_A)$ be the bounded derived category $\D^b \mod (\CC Q_A)$ of finitely generated modules over the path algebra $\CC Q_A$.
\begin{prop}[{\cite[Proposition 2.4]{GL}, cf.~\cite[Proposition 2.24]{STW}}]\label{prop : derived equivalence}
There exist triangle equivalences $\D_A \cong \D^b (Q_A)$. 
\qed
\end{prop}

We consider a generalized root system associated with the derived category $\D_A$.
\begin{prop}\label{prop : GRS from derived categories}
Let $\E = (E_1, \dots, E_{\mu_A})$ be a full exceptional collection in $\D_A$.
The tuple $(R_{\D_A}, \mathbf{c}_{\D_A}) = (L, I, \Delta_\mathrm{re}, \mathbf{c})$ consisting of 
\begin{itemize}
\item the Grothendieck group $L \coloneqq K_0(\D_A)$,
\item $I \coloneqq \chi_{\D_A} + \chi^T_{\D_A}$, where $\chi_{\D_A}$ is the Euler form on $K_0(\D_A)$ defined by 
\begin{equation*}
\chi_{\D_A}([X], [Y]) \coloneqq \sum_{p \in \ZZ} (-1)^p \dim_\CC \Hom_{\D_A} (X, Y[p])
\end{equation*}
for $X, Y \in \D_A$,
\item the subset $\Delta_\mathrm{re}$ of $L$ defined by 
\begin{equation*}
\Delta_\mathrm{re} \coloneqq \{ w([E_i]) \in L \mid w \in W(\E), \ i = 1, \dots, {\mu_A} \},
\end{equation*}
where $W(\E)$ is the subgroup of $\Aut_\ZZ (L, I)$ generated by reflections $r_{[E_i]}$ of $[E_i]$,
\item the automorphism $\mathbf{c}$ on $L$ induced by the $(-1)$-shifted Serre functor $\S_{\D_A}[-1]$,
\end{itemize}
is isomorphic to a generalized root system $(\widetilde{R}_A, \widetilde{\mathbf{c}}_A)$ of the corresponding affine ADE type.

Moreover, the generalized root system $(R_{\D_A}, \mathbf{c}_{\D_A})$ does not depend on the choice of a full exceptional collection.
In particular, the Coxeter transformation $\mathbf{c}_{\D_A}$ satisfies 
\begin{equation*}
\mathbf{c}_{\D_A} = r_{[E_1]} \cdots r_{[E_{\mu_A}]},
\end{equation*}
for any full exceptional collection $(E_1, \dots, E_{\mu_A})$.
\end{prop}
\begin{pf}
By Proposition \ref{prop : derived equivalence} and \cite{C-B, R}, the triangulated category $\D_A$ satisfies the assumptions of \cite[Proposition 2.10, Lemma 2.11]{STW}.
\qed
\end{pf}

We identify the generalized root system $(R_{\D_A}, \mathbf{c}_{\D_A})$ with $(\widetilde{R}_A, \widetilde{\mathbf{c}}_A)$.
The Euler form $\chi_{\D_A} \colon K_0(\D_A) \times K_0(\D_A) \longrightarrow \ZZ$ coincides with the one $\chi \colon \widetilde{L} \times \widetilde{L} \longrightarrow \ZZ$ of $(\widetilde{R}_A, \widetilde{\mathbf{c}}_A)$.
\begin{rem}
The generalized root system obtained from a triangulated category admits a standard Euler form.
One can define twist automorphisms on the generalized root system by this standard Euler form.
On the other hand, the uniqueness of the Euler form in the sense of Definition \ref{defn : Euler form} does not hold in general.
Nevertheless, the generalized root systems of affine ADE type admit a unique Euler form (see Lemma \ref{lem : uniqueness of the Euler from}).
\end{rem}

Let $S_v$ denote the simple $\CC \widetilde{\TT}_A$-module corresponding to the vertex $v \in \widetilde{V}_A$.
A direct calculation shows that the ordered set $(S_\mathbf{1}, S_{(1,1)}, \dots, S_{(1, a_1 - 1)}, S_{(2, 1)}, \dots, S_{(2, a_2 - 1)}, S_{(3, 1)},\dots, S_{(3, a_3 - 1)}, S_{\mathbf{1}^\ast})$ forms a full exceptional collection. 
A real root $\alpha \in \widetilde{\Delta}_\mathrm{re}$ is said to be {\em positive} if $\alpha = \sum_{v \in \widetilde{V}_A} c_v [S_v]$ satisfies $c_v \ge 0$ for all $v \in \widetilde{V}_A$.
Denote by $\widetilde{\Delta}^+_\mathrm{re}$ the set of positive real roots.
Proposition \ref{prop : derived equivalence} implies the following 
\begin{prop}[{\cite[Theorem 1]{K1}}]\label{prop : Kac theorem}
For any $\alpha \in \widetilde{\Delta}^+_\mathrm{re}$, there exists exactly one indecomposable object $E \in \D_A$ up to shift and isomorphism such that $[E] = \alpha$ on $K_0(\D_A)$.
\qed
\end{prop}

\begin{prop}\label{prop : spherical twist and twist automorphism}
The following holds:
\begin{enumerate}
\item An object $S \in \D_A$ is $1$-spherical if and only if $[S] = \delta$ or $-\delta$.
\item For a spherical object $S \in \D_A$, we have 
\begin{equation}\label{eq : spherical twist and twist auto}
[\T_S (X)] = t_\delta ([X]), 
\end{equation}
for any $X \in \D_A$.
\item We have $t_\delta^{- \chi_A \ell_A} = \widetilde{\mathbf{c}}_A^{\ell_A}$, where $\ell_A = \mathrm{lcm}(a_1, a_2, a_3)$.
\end{enumerate}
\end{prop}
\begin{pf}
(1): It follows from \cite[Section 6]{DR} and \cite[Section 5.5]{Ha} that a module $S \in \mod(\CC Q_A)$ is $1$-spherical if and only if $[S] = \delta$.
Since the algebra $\CC Q_A$ is hereditary, any $1$-spherical object in $\D^b(Q_A)$ is isomorphic to a shifted $1$-spherical object in $\mod(\CC Q_A)$.
Hence, Proposition \ref{prop : derived equivalence} implies the statement.

(2): By the exact triangle \eqref{eq : spherical twist}, we have 
\begin{eqnarray*}
[\T_S (X)] & = & [X] - \big[ \Hom^\bullet_{\D_A} (S, X) \otimes S \big] \\
& = & [X] - \Big[ \bigoplus_{p \in \ZZ} S[p]^{\oplus \dim_\CC \Hom (S, X[p])} \Big] \\
& = & [X] - \sum_{p \in \ZZ} (- 1)^p \dim_\CC \Hom_{\D_A} (S, X[p]) [S] \\
& = & [X] - \chi([S], [X]) [S].
\end{eqnarray*}
Hence, the statement (1) yields equality \eqref{eq : spherical twist and twist auto}.

(3): Combining \cite[Proposition 4.8]{OST} with \cite[Section 2.2]{GL}, we obtain $\T_S^{ -\chi_A \ell_A} \cong \S_{\D_A}^{\ell_A} [- \ell_A]$ for any $1$-spherical object $S \in \D_A$.
Thus, the statement follows from (2).
\qed
\end{pf}
\subsection{Lyashko--Looijenga map and root bases}
In this section, we show that the degree of the Lyashko--Looijenga map is equal to the number of root bases which give the Coxeter transformation $\widetilde{\mathbf{c}}_A$ modulo the twist automorphism.

The {\em Artin's braid group} $\Br_{\mu_A}$ on $\mu_A$-stands is a group presented by the following generators and relations: 
\begin{description}
\item[{\bf Generators}] $\{ \sigma_i \mid i=1, \dots, \mu_A - 1 \}$
\item[{\bf Relations}] $\sigma_{i} \sigma_{j} = \sigma_{j} \sigma_{i}$ for $|i - j| \ge 2$, $\sigma_{i} \sigma_{i+1} \sigma_{i} = \sigma_{i+1} \sigma_{i} \sigma_{i+1}$ for $i = 1, \dots, \mu_A - 2$.
\end{description}
Consider the semi-direct product $\Br_ {\mu_A} \ltimes \ZZ^{\mu_A}$, which is defined by $\Br_ {\mu_A} \rightarrow \mathfrak{S}_{\mu_A} \rightarrow \Aut_\ZZ \ZZ^{\mu_A}$, where the first homomorphism is $b_i \mapsto (i, i+1)$ and the second one is induced by the natural actions of the symmetric group $\mathfrak{S}_{\mu_A}$ on $\ZZ^{\mu_A}$. 
The group $\Br_{\mu_A} \ltimes \ZZ^{\mu_A}$ acts on $\FEC(\D_A)$ by 
\begin{gather*}
\sigma_{i} \cdot (E_{1}, \dots, E_{{\mu_A}}) \coloneqq (E_{1}, \dots, E_{i-1}, E_{i+1}, \RR_{E_{i+1}} E_i, E_{i+2}, \dots, E_{{\mu_A}}), \\
\sigma^{-1}_{i} \cdot (E_{1}, \dots, E_{{\mu_A}}) \coloneqq (E_{1}, \dots, E_{i-1}, \LL_{E_{i}} E_{i+1}, E_{i}, E_{i+2}, \dots, E_{{\mu_A}}), \\
e_{i} \cdot (E_{1}, \dots, E_{{\mu_A}}) \coloneqq (E_{1}, \dots, E_{i-1}, E_{i} [1], E_{i+1}, \dots, E_{{\mu_A}}),
\end{gather*}
where $\RR_F E$ (resp. $\LL_E F$) is called the {\em right mutation} of $E$ through $F$ (resp. {\em left mutation} of $F$ through $E$) defined by 
\begin{subequations}
\begin{equation}\label{eq : right mutation}
\RR_F E \longrightarrow E \overset{{\rm ev}^*}{\longrightarrow} \RR \Hom_{\D_A}(E, F)^* \otimes F,
\end{equation}
\begin{equation}\label{eq : left mutation}
\RR\Hom_{\D_A}(E, F) \otimes E \overset{\rm ev}{\longrightarrow} F \longrightarrow \LL_E F,
\end{equation}
\end{subequations}
and $e_{i}$ is the $i$-th generator of $\ZZ^{\mu_A}$.
\begin{prop}[{\cite[Theorem]{C-B}}]\label{prop : transitivity by C-B}
The $\Br_ {\mu_A} \ltimes \ZZ^{\mu_A}$-action on $\FEC(\D_A)$ is transitive.
\qed
\end{prop}

For the generalized root system $(\widetilde{R}_A, \widetilde{\mathbf{c}}_A)$ of affine ADE type, we define the set $\B_A$ by 
\begin{equation*}
\B_A \coloneqq \{ (\alpha_1, \dots, \alpha_{\mu_A}) \mid \widetilde{\mathbf{c}}_A = r_{\alpha_1} \cdots r_{\alpha_{\mu_A}}, \ \text{$B = \{ \alpha_1, \dots, \alpha_{\mu_A} \}$ is a root basis} \}.
\end{equation*}
Then, there is a group action of $\Br_{\mu_A} \ltimes (\ZZ/2\ZZ)^{\mu_A}$ on $\B_A$ given by 
\begin{gather*}
\sigma_i (\alpha_1, \dots, \alpha_{\mu_A}) \coloneqq (\alpha_1, \dots, \alpha_{i - 1}, \alpha_{i + 1}, r_{\alpha_{i + 1}} \alpha_i, \alpha_{i + 2}, \dots, \alpha_{\mu_A}),\\
\sigma^{-1}_i (\alpha_1, \dots, \alpha_{\mu_A}) \coloneqq (\alpha_1, \dots, \alpha_{i - 1}, r_{\alpha_i} \alpha_{i + 1}, \alpha_i, \alpha_{i + 2}, \dots, \alpha_{\mu_A}),\\
e_i (\alpha_1, \dots, \alpha_{\mu_A}) \coloneqq (\alpha_1, \dots, \alpha_{i - 1}, - \alpha_i, \alpha_{i + 1}, \dots, \alpha_{\mu_A}).
\end{gather*}
\begin{prop}[{\cite[Theorem 1.4]{IS}}]\label{prop : transitivity by IS}
The $\Br_{\mu_A} \ltimes (\ZZ/2\ZZ)^{\mu_A}$-action on $\B_A$ is transitive.
\qed
\end{prop}

By Proposition \ref{prop : GRS from derived categories}, one can define a map from the set of full exceptional collections to the set of root basis whose products of reflections is the Coxeter transformation.
\begin{prop}[{\cite[Proposition 4.6]{HR}, cf.~\cite{IS}, \cite[Corollary 4.9]{BWY}}]\label{prop : one-to-one correspondence}
The map 
\begin{equation}\label{eq : one-to-one correspondence}
\FEC(\D_A) / \ZZ^{\mu_A} \longrightarrow \B_A / (\ZZ / 2 \ZZ)^{\mu_A}, \quad [(E_1, \dots, E_{\mu_A})] \mapsto [(\varepsilon_1, \dots, \varepsilon_{\mu_A})]
\end{equation}
is a $\Br_{\mu}$-equivariant bijection, where $\varepsilon_i \coloneqq [E_i] \in L$ for $i = 1, \dots, {\mu_A}$.
\qed
\end{prop}

Define a full exceptional collection $\E_A = (E_{\mathbf{1}^\ast}, E_\mathbf{1}, E_{(1,1)}, \dots, E_{(3, a_3 - 1)})$ by 
\begin{equation*}
\E_A \coloneqq \sigma_1^{-1} \sigma_2^{-1} \cdots \sigma_{\mu_A - 1}^{-1} \cdot (S_\mathbf{1}, S_{(1,1)}, \dots, S_{(3, a_3 - 1)}, S_{\mathbf{1}^\ast}).
\end{equation*}
\begin{cor}\label{cor : our full exceptional collection}
The map \eqref{eq : one-to-one correspondence} sends $\E_A$ to $\widetilde{B}_A$.
\end{cor}
\begin{pf}
By the definition of the left mutation, a direct calculation yields the statement.
\qed
\end{pf}

Denote by $\ST_1(\D_A) \coloneqq \br{\T_S \in \Aut(\D_A) \mid \text{$S$ is $1$-spherical} }$ the subgroup generated by $1$-spherical twists on $\D_A$.
There is an action on $\FEC(\D_A)$ defined by 
\begin{equation*}
\T_S \cdot (E_{1}, \dots, E_{\mu_A}) \coloneqq (\T_S(E_{1}), \dots, \T_S(E_{\mu_A})).
\end{equation*}
According to Proposition \ref{prop : transitivity by C-B}, one can consider the subgroup $\br{\ZZ^{\mu_A}, \mathrm{ST}(\D_A)}$ of $\Br_{\mu_A} \ltimes \ZZ^{\mu_A}$ generated by $\ZZ^{\mu_A}$ and $\ST_1(\D_A)$.
Similarly, we have a $T_A$-action on $\B_A$ defined by 
\begin{equation*}
t_\delta \cdot (\alpha_1, \dots, \alpha_{\mu_A}) \coloneqq (t_\delta(\alpha_1), \dots, t_\delta(\alpha_{\mu_A})).
\end{equation*}
By Proposition \ref{prop : transitivity by IS}, one can also consider the subgroup $\br{(\ZZ / 2\ZZ)^{\mu_A}, T_A}$ of $\Br_{\mu_A} \ltimes (\ZZ / 2\ZZ)^{\mu_A}$ generated by $(\ZZ / 2\ZZ)^{\mu_A}$ and $T_A$.
It was proved by \cite[Proposition 4.8]{OST} that $\ST_1(\D_A) \cong \ZZ$.
Thus, the following corollary follows from Proposition \ref{prop : spherical twist and twist automorphism} and Proposition \ref{prop : one-to-one correspondence}.
\begin{cor}\label{cor : main corollary}
The map 
\begin{equation*}
\FEC(\D_A) / \br{\ZZ^{\mu_A}, \ST_1(\D_A)} \longrightarrow \B_A / \br{\ZZ_2^{\mu_A}, T_A}, \quad [(E_1, \dots, E_{\mu_A})] \mapsto [(\varepsilon_1, \dots, \varepsilon_{\mu_A})]
\end{equation*}
is a bijection.
\qed
\end{cor}

Let $e(\B_A)$ be the cardinality of $\B_A / \br{\ZZ_2^\mu, T_A}$ and $e(\D_A)$ the one of $\FEC(\D_A) / \br{\ZZ^\mu, \ST_1(\D_A)}$.
It was proved by \cite[Theorem 1.2]{OST} that $\deg LL = e(\D_A)$.
Since Corollary \ref{cor : main corollary} implies $e(\D_A) = e (\B_A)$, we obtain the following
\begin{cor}\label{cor : LL map and root bases}
We have
\begin{equation*}
\deg LL = e(\B_A) = e(\D_A).
\end{equation*}
\qed
\end{cor}
This corollary can be regarded as an analogue of Deligne's result \cite{De} for finite ADE cases.

\section{Extended Artin group and spherical twists}\label{sec : Extended Artin group and spherical twists}
In this section, we introduce the extended Artin group associated with the Coxeter--Dynkin diagram $\widetilde{\Gamma}_A$ (see Figure \ref{fig : Coxeter Dynkin diagram associated with the octopus quiver}).
We prove that the group is isomorphic to the fundamental group of the regular subspace of the Frobenius manifold in Corollary \ref{cor : Frobenius structure}.
We conjecture the relation between the monodromy representation of the Frobenius manifold and a representation of the group generated by spherical twists on an $\XX$-Calabi--Yau category. 
\subsection{Extended Artin group}
In order to introduce an extended Artin group associated with $(\widetilde{R}_A, \widetilde{\mathbf{c}}_A)$, we consider the following Artin group defined by Shiraishi--Takahashi--Wada.
\begin{defn}[{\cite[Definition 5.1]{STW}}]
Define the {\em Artin group} $\widetilde{G}_A$ associated with the Coxeter--Dynkin diagram $\widetilde{\Gamma}_A$ as follows: 
\begin{itemize}
\item {\bf Generators:} $\{ \sigma_v \mid v \in \widetilde{V}_A \}$, 
\item {\bf Relations:} 
\begin{subequations}\label{eq : affine Artin relation}
\begin{alignat}{2}
\sigma_v \sigma_{v'} & = \sigma_{v'} \sigma_v, & \quad \quad & \nexists \ \text{edge between $v$ and $v'$} \\
\sigma_v \sigma_{v'} \sigma_v & = \sigma_{v'} \sigma_v \sigma_{v'}, & & \exists \ \text{edge between $v$ and $v'$} \\
\rho_\mathbf{1} \rho_{(i, 1)} & = \rho_{(i, 1)} \rho_\mathbf{1}, & & i = 1, 2, 3, \\
\sigma_{(i, 1)} \rho_{(i, 1)} & = \rho_{(i, 1)} \sigma_{(i,1)}, & & i = 1, 2, 3,
\end{alignat}
\end{subequations}
where we put $\rho_\mathbf{1} = \sigma_\mathbf{1} \sigma_{\mathbf{1}^\ast}$ and $\rho_{(i, 1)} = \sigma_{(i,1)} \rho_\mathbf{1} \sigma_{(i,1)} \rho_\mathbf{1}^{-1}$.
\end{itemize}
\end{defn}
On the other hand, one can define the affine Artin group for the affine Weyl group.
For a root system $\widetilde{R}_A$ of affine ADE type, denote by $G(\widetilde{R}_A)$ the affine Artin group (see \cite{BS}).
\begin{prop}[{\cite[Theorem 5.8]{STW}}]\label{prop : affine Artin group}
There exists an isomorphism $\widetilde{G}_A \cong G(\widetilde{R}_A)$.
\qed
\end{prop}

\begin{cor}\label{cor : fundamental group and Artin group}
We have $\pi_1(\X_A^\mathrm{reg} / \widetilde{W}_A) \cong \widetilde{G}_A$.
\end{cor}
\begin{pf}
By \cite[Theorem 2.20]{I}, we have $\pi_1(X_A^\mathrm{reg} / \widetilde{W}_A) \cong \pi_1(X_A) \times G(\widetilde{R}_A)$.
Hence, the statement follows from Proposition \ref{prop : affine Artin group}.
\qed
\end{pf}

Based on Lemma \ref{lem : basic properties}, we introduce an extended Artin group associated.
\begin{defn}\label{defn : extended Artin group}
Define the {\em extended Artin group} $\widehat{G}_A$ associated with the Coxeter--Dynkin diagram $\widetilde{\Gamma}_A$ as follows: 
\begin{itemize}
\item {\bf Generators:} $\{ \sigma_v \mid v \in \widetilde{\Gamma}_A \} \cup \{ \tau \}$, 
\item {\bf Relations:} \eqref{eq : affine Artin relation} and 
\begin{subequations}
\begin{alignat}{2}
\label{eq : extended Artin 1} \tau \sigma_{\mathbf{1}^\ast} \tau^{-1} & = \sigma_{\mathbf{1}}, & \\
\label{eq : extended Artin 2} \tau \sigma_{\mathbf{1}} \tau^{-1} & = \sigma^{-1}_{\mathbf{1}} \sigma_{\mathbf{1}^\ast} \sigma_{\mathbf{1}}, \quad & \\
\label{eq : extended Artin 3} \tau \sigma_v \tau^{-1} & = \sigma_v, \quad & v \in \widetilde{V}_A \setminus \{ \mathbf{1}^\ast, \mathbf{1} \}.
\end{alignat}
\end{subequations}
\end{itemize}
\end{defn}

\begin{rem}
It is known that a root system $\widetilde{R}$ of affine $A_\mu$ type is obtained from an annulus with $(\mu + 1)$-decoration points.
The {\em annular braid group} $CB_{\mu + 1}$ is defined to be the fundamental group of the configuration space of $\mu + 1$ distinct points in $\CC^\ast$ (cf.~\cite{KP}).
For a tuple $A = (1, p, q)$ with $p + q = \mu$, one can check that $\widetilde{G}(\widetilde{R}) \subsetneq \widehat{G}_A \subsetneq CB_{\mu + 1}$ (cf.~Section \ref{sec : Spherical twists group}).
\end{rem}

\begin{prop}
The map $\widehat{G}_A \longrightarrow \widehat{W}_A$ given by $\sigma_v \mapsto r_{\alpha_v}$ and $\tau \mapsto (t_\delta, \gamma)$ is surjective.
\end{prop}
\begin{pf}
It was proved by \cite[Proposition 5.4]{STW} that the map $\widetilde{G}_A \longrightarrow \widetilde{W}_A$ defined by $\sigma_v \mapsto r_{\alpha_v}$ is a surjective homomorphism.
By Lemma \ref{lem : basic properties}, there is a following commutative diagram
\begin{equation*}
\xymatrix{
1 \ar[r] & \widetilde{G}_A \ar[r] \ar@{->>}[d] & \widehat{G}_A \ar[r] \ar[d] & \br{\tau} \ar[r] \ar[d]^{\cong} & 1 \\
1 \ar[r] & \widetilde{W}_A \ar[r] & \widetilde{W}_A^\mathrm{ext} \ar[r] & T_A \ar[r] & 1
}
\end{equation*}
where the map $\widehat{G}_A \longrightarrow \widetilde{W}_A^\mathrm{ext}$ is defined by $\sigma_v \mapsto r_{\alpha_v}$ and $\tau \mapsto t_\delta$.
It follows from the five lemma that the map is surjective.
Thus, we obtain the statement.
\qed
\end{pf}

\begin{prop}\label{prop : fundamental group and extended Artin group}
We have $\pi_1(\X_A^\mathrm{reg} / \widehat{W}_A) \cong \widehat{G}_A$.
\end{prop}
\begin{pf}
Since the covering map $\X_A^\mathrm{reg} / \widetilde{W}_A \longrightarrow \X_A^\mathrm{reg} / \widehat{W}_A$ is regular, we have a commutative diagram
\begin{equation*}
\xymatrix{
1 \ar[r] & \pi_1(\X_A^\mathrm{reg} / \widetilde{W}_A) \ar[r] \ar@{->>}[d] & \pi_1(\X_A^\mathrm{reg} / \widehat{W}_A) \ar[r] \ar@{->>}[d] & C_A \ar[r] \ar@{=}[d] & 1 \\
1 \ar[r] & \widetilde{W}_A \ar[r] & \widehat{W}_A \ar[r] & C_A \ar[r] & 1
}
\end{equation*}
Hence, by Corollary \ref{cor : fundamental group and Artin group}, there is a map $\pi_1(\X_A^\mathrm{reg} / \widehat{W}_A) \longrightarrow \widehat{G}_A$ such that the following diagram commutes:
\begin{equation*}
\xymatrix{
1 \ar[r] & \pi_1(\X_A^\mathrm{reg} / \widetilde{W}_A) \ar[r] \ar[d]^{\cong} & \pi_1(\X_A^\mathrm{reg} / \widehat{W}_A) \ar[r] \ar[d] & C_A \ar[r] \ar[d]^{\cong} & 1 \\
1 \ar[r] & \widetilde{G}_A \ar[r] & \widehat{G}_A \ar[r] & \br{\tau} \ar[r] & 1
}
\end{equation*}
Thus, the five lemma yields the statement.
\qed
\end{pf}

As a direct consequence, we can describe the map \eqref{defn : monodromy of Frobenius manifold} by the extended Artin group.
\begin{cor}\label{cor : extended Artin group and extended Weyl group}
There is a surjective homomorphism $\widehat{G}_A \longrightarrow \widehat{W}_A$ such that the following diagram commutes:
\begin{equation*}
\xymatrix{
\widehat{G}_A \ar@{->>}[d] \ar[rr]^-\cong & & \pi_1(\X_A^\mathrm{reg} / \widehat{W}_A) \ar[d]^{\rho} \\
\widehat{W}_A \ar[rr] & & \Aut(\X_A) 
}
\end{equation*}
\qed
\end{cor}
\subsection{Extended Seidel--Thomas braid group}\label{sec : Spherical twists group}
We consider an $\XX$-Calabi--Yau triangulated category associated with the algebra $\CC \widetilde{\TT}_A$.
We refer to \cite{Ke, IQ} for the Calabi--Yau completion.
Let $\Theta_{\CC \widetilde{\TT}_A}$ be the inverse dualizing complex of $\CC \widetilde{\TT}_A$.
The {\em $\XX$-Calabi--Yau completion} of $\CC \widetilde{\TT}_A$ is the differential $(\ZZ \oplus \ZZ \XX)$-graded $\CC$-algebra $\Pi_\XX (\CC \widetilde{\TT}_A)$ defined by 
\begin{equation*}
\Pi_\XX (\CC \widetilde{\TT}_A) \coloneqq (\CC \widetilde{\TT}_A) \oplus \theta \oplus ( \theta \otimes_{\CC \widetilde{\TT}_A}  \theta) \oplus \cdots,
\end{equation*}
where $\theta = \Theta_{\CC \widetilde{\TT}_A}[\XX - 1]$.
For simplicity, denote by $\D^\XX_A$ the full subcategory of the derived category $\D(\Pi_\XX (\CC \widetilde{\TT}_A))$ consisting of objects with finite-dimensional total cohomology.
It follows from \cite[Theorem 2.4]{IQ} (cf.~\cite{Ke}) that $\D^\XX_A$ is $\XX$-Calabi--Yau, namely the Serre functor is naturally isomorphic to $[\XX]$.
The following lemma is useful to study the relation between a spherical twist and an autoequivalence:
\begin{lem}[{cf.~\cite[Lemma 8.21]{Huy}}]\label{lem : autoequivalence and spherical twist}
Let $S^\XX \in \D^\XX_A$ be an $\XX$-spherical object.
For any autoequivalence $\Phi \in \Aut (\D^\XX_A)$, we have $\Phi \circ \T_{S^\XX} \circ \Phi^{-1} \cong \T_{\Phi(S^\XX)}$.
\qed
\end{lem}

The canonical projection $\Pi_\XX(\CC \widetilde{\TT}_A) \longrightarrow \CC \widetilde{\TT}_A$ induces an exact functor $\L \colon \D_A \longrightarrow \D^\XX_A$ such that for $X, Y \in \D^\XX_A$ there is a canonical isomorphism 
\begin{equation}\label{eq : Lagrangian immersion}
\Hom^\bullet_{\D^\XX_A} (\L(X), \L(Y)) \cong \Hom^\bullet_{\D_A} (X, Y) \oplus \Hom^\bullet_{\D_A} (Y, X)^\ast [-\XX].
\end{equation}
The functor $\L \colon \D_A \longrightarrow \D^\XX_A$ is called the {\em Lagrangian immersion} \cite{IQ, KQ}.

\begin{lem}\label{lem : mutation and spherical twist}
Let $(E, F)$ be an exceptional pair in $\D_A$ and put $E^\XX = \L(E), F^\XX = \L(F) \in \D^\XX_A$.
We have $\L(\RR_E F) \cong \T^{-1}_{E^\XX} (F^\XX)$ and $\L(\LL_F E) \cong \T_{F^\XX}(E^\XX)$.
\end{lem}
\begin{pf}
By equation \eqref{eq : Lagrangian immersion}, we have $\Hom^\bullet_{\D^\XX_A} (E^\XX, F^\XX) \cong \Hom^\bullet_{\D_A} (E, F)$.
There is a commutative diagram 
\begin{equation*}
\xymatrix{
\L(E) \ar[rr] \ar@{=}[d] & & \L \big( \Hom^\bullet_{\D_A} (E, F)^\ast \otimes F \big) \ar[d]^{\cong} \\
E^\XX \ar[rr] & & \Hom^\bullet_{\D^\XX_A} (E^\XX, F^\XX)^\ast \otimes F^\XX,
}
\end{equation*}
which yields $\L(\RR_E F) \cong \T^{-1}_{E^\XX} (F^\XX)$.

One can also prove $\L(\LL_F E) \cong \T_{F^\XX}(E^\XX)$ in the similar way.
\qed
\end{pf}

\begin{lem}\label{lem : Lagrangian immersion}
For each $\Phi \in \Aut(\D_A)$, there exists an autoequivalence $\Phi^\L \in \Aut(\D^\XX_A)$ such that $\L \circ \Phi = \Phi^\L \circ \L$.
Especially, there is an injection $(-)^\L \colon \Aut(\D_A) \longrightarrow \Aut(\D^\XX_A)$.
\end{lem}
\begin{pf}
The proof is done in the same way as in \cite[Lemma A.8]{Q}.
\qed
\end{pf}

Let $E_v \in \D_A$ be an exceptional object given in Corollary \ref{cor : our full exceptional collection} and $E_v^\XX \coloneqq \L (E_v)$.
By \eqref{eq : Lagrangian immersion}, the object $E_v^\XX$ is an $\XX$-spherical object in $\D^\XX_A$.
The {\em Seidel--Thomas braid group} $\ST^\circ(\D^\XX_A)$ of the triangulated category $\D^\XX_A$ (with respect to $\widetilde{V}_A$) is defined to be the subgroup of $\Aut (\D^\XX_A)$ generated by the spherical twists $\T_{E_v^\XX}$:
\begin{equation*}
\ST^\circ(\D^\XX_A) \coloneqq \br{\T_{E^\XX_v} \mid v \in \widetilde{V}_A}
\end{equation*}
For the extended Dynkin quiver $Q_A$, one can define the $\XX$-Calabi--Yau triangulated category $\D^\XX(Q_A)$ and $\ST^\circ(\D^\XX (Q_A))$ in the same way.

\begin{lem}[{cf.~\cite[Theorem 6.12]{STW}}]\label{lem : artin group and Seidel--Thomas braid group}
The map $\iota_\XX \colon \widetilde{G}_A \longrightarrow \ST^\circ(\D^\XX_A)$ is a surjective group homomorphism.
\end{lem}
\begin{pf}
We first show $\ST^\circ(\D^\XX_A) = \ST^\circ(\D^\XX(Q_A))$.
Let $F_i$ denote the simple $\CC Q_A$-module corresponding to the vertex $i$.
Since the ordered set $(F_1, \dots, F_{\mu_A})$ forms a full exceptional collection, Proposition \ref{prop : transitivity by C-B} implies that there exists an element $b \in \Br_{\mu_A} \ltimes \ZZ^{\mu_A}$ such that $\E_A = b \cdot (T_1, \dots, T_{\mu_A})$.
Especially, each $E_v$ is obtained by iteration of mutations and shifts from some $F_i$.
Hence, it follows from Proposition \ref{prop : spherical object and twist}, Lemma \ref{lem : mutation and spherical twist} and Proposition \ref{lem : autoequivalence and spherical twist} that $\T_{E_v^\XX} \in \ST^\circ(\D^\XX (Q_A))$, and we have $\ST^\circ(\D^\XX_A) \subset \ST^\circ(\D^\XX(Q_A))$.
similarly, we have $\ST^\circ(\D^\XX(Q_A)) \subset \ST^\circ(\D^\XX_A)$.

By \cite[Section A.4]{IQ}, the map $G(\widetilde{R}_A) \longrightarrow \ST^\circ(\D^\XX (Q_A))$ defined by \cite[Equation (A.1)]{IQ} is surjective.
Since there is a commutative diagram 
\begin{equation*}
\xymatrix{
\widetilde{G}_A \ar[d]_{\iota_\XX} \ar[rr]^{\cong} & & G(\widetilde{R}_A) \ar@{->>}[d] \\
\ST^\circ(\D^\XX_A) \ar@{=}[rr] & &  \ST^\circ(\D^\XX(Q_A)),
}
\end{equation*}
we obtain the statement.
\qed
\end{pf}

As an analogue of Definition \ref{defn : extended affine Weyl group},  we introduce the following 
\begin{defn}\label{defn : extended Seidel--Thomas braid group}
Define a subgroup $\widehat{\ST}(\D^\XX_A)$ of $\Aut (\D^\XX_A)$ by 
\begin{equation}
\widehat{\ST}(\D^\XX_A) \coloneqq \big\langle \ST^\circ(\D^\XX_A), \ \ST^\L_1(\D_A) \big\rangle,
\end{equation}
where $\ST^\L_1(\D_A) \coloneqq \br{\T^\L_S \in \Aut(\D^\XX_A) \mid \text{$S \in \D_A$ is $1$-spherical}}$.
We call $\widehat{\ST}(\D^\XX_A)$ the {\em extended Seidel--Thomas braid group} of $\D^\XX_A$.
\end{defn}

\begin{prop}\label{prop : relation of the extended part}
Let $S \in \D_A$ be a $1$-spherical object.
We have 
\begin{subequations}
\begin{alignat}{2}
\label{eq : relation of extended braid group 1} \T^\L_S \circ \T_{E_{\mathbf{1}^\ast}^\XX} \circ (\T^\L_S)^{-1} & \cong  \T_{E^\XX_{\mathbf{1}}}, & \\
\label{eq : relation of extended braid group 2} \T^\L_S \circ \T_{E_{\mathbf{1}}^\XX} \circ (\T^\L_S)^{-1} & \cong \T^{-1}_{E^\XX_{\mathbf{1}}} \circ \T_{E^\XX_{\mathbf{1}^\ast}} \circ \T_{E^\XX_{\mathbf{1}}}, \quad & \\
\label{eq : relation of extended braid group 3} \T^\L_S \circ \T_{E_v^\XX} \circ (\T^\L_S)^{-1} & \cong \T_{E^\XX_v}, \quad & v \in \widetilde{V}_A \setminus \{ \mathbf{1}^\ast, \mathbf{1} \}. 
\end{alignat}
\end{subequations}
\end{prop}
\begin{pf}

We first show relation \eqref{eq : relation of extended braid group 1}.
By Proposition \ref{prop : spherical object and twist} and Lemma \ref{lem : autoequivalence and spherical twist}, it is sufficient to show $\T^\L_S (E^\XX_{\mathbf{1}^\ast}) \cong E^\XX_{\mathbf{1}}[n]$ for some $n \in \ZZ$. 
Lemma \ref{lem : basic properties} (4), Proposition \ref{prop : spherical twist and twist automorphism} and Corollary \ref{cor : our full exceptional collection} yield that 
\begin{equation*}
[\T_S(E_{\mathbf{1}^\ast})] = t_\delta (\alpha_{\mathbf{1}^\ast}) = \alpha_{\mathbf{1}} = [E_\mathbf{1}] \quad \text{on $K_0(\D_A)$}.
\end{equation*}
It follows from Proposition \ref{prop : Kac theorem} that there is an integer $n \in \ZZ$ such that $\T_S(E_{\mathbf{1}^\ast}) \cong E_{\mathbf{1}} [n]$.
Thus, we have $\T^\L_S (E^\XX_{\mathbf{1}^\ast}) \cong \L(\T_S(E_{\mathbf{1}^\ast})) \cong E^\XX_{\mathbf{1}} [n]$.

Next, we prove $\T^\L_S (E^\XX_{\mathbf{1}}) \cong \T^{-1}_{E^\XX_{\mathbf{1}}} (E^\XX_{\mathbf{1}^\ast}) [n]$ for some $n \in \ZZ$ to see relation \eqref{eq : relation of extended braid group 2}.
Since $[\RR_{E_\mathbf{1}} E_{\mathbf{1}^\ast}] = r_{\alpha_\mathbf{1}} (\alpha_{\mathbf{1}^\ast})$, we have $[\T_S(E_{\mathbf{1}})] = [\RR_{E_\mathbf{1}} E_{\mathbf{1}^\ast}]$.
It follows again from Proposition \ref{prop : Kac theorem} that there is an integer $n \in \ZZ$ such that $\T_S(E_{\mathbf{1}}) \cong \RR_{E_\mathbf{1}} E_{\mathbf{1}^\ast} [n]$.
By Lemma \ref{lem : mutation and spherical twist}, we obtain $\T^\L_S (E^\XX_{\mathbf{1}}) \cong \L(\RR_{E_\mathbf{1}} E_{\mathbf{1}^\ast} [n]) \cong \T^{-1}_{E^\XX_{\mathbf{1}}} (E^\XX_{\mathbf{1}^\ast}) [n]$.

Similarly, one can check easily that $\T^\L_S(E^\XX_{(i,j)}) \cong E^\XX_{(i,j)}$ for $i = 1, 2, 3$ and $j = 1, \dots, a_i - 1$.
\qed
\end{pf}

\begin{thm}\label{thm : extended Artin group and extended ST group}
The map 
\begin{equation*}
\widehat{\iota}_\XX \colon \widehat{G}_A \longrightarrow \widehat{\ST}(\D^\XX_A), \quad 
\begin{cases}
\sigma_v \mapsto\T_{E^\XX_v}, \quad v \in \widetilde{V}_A, \\
\tau \mapsto \T^\L_S,
\end{cases}
\end{equation*}
is a surjective homomorphism.
\end{thm}
\begin{pf}
Lemma \ref{lem : artin group and Seidel--Thomas braid group} and Proposition \ref{prop : relation of the extended part} implies that there is a commutative diagram whose rows are short exact sequences
\begin{equation}\label{eq : Artin and braid}
\xymatrix{
1 \ar[r] & \widetilde{G}_A \ar[r] \ar@{->>}[d]^{\iota_\XX} & \widehat{G}_A \ar[r] \ar[d]^{\widehat{\iota}_\XX} & \br{\tau} \ar[r] \ar[d]^\cong & 1 \\ 
1 \ar[r] & \ST^\circ(\D^\XX_A) \ar[r] & \widehat{\ST}(\D^\XX_A) \ar[r] & \ST_1^\L(\D_A) \ar[r] & 1 
}.
\end{equation}
The five lemma yields the statement.
\qed
\end{pf}

Therefore, we obtain a commutative diagram
\begin{equation*}
\xymatrix{
\pi_1(\X_A^\mathrm{reg} / \widehat{W}_A) \ar@{->>}[rr] \ar[d]_\rho & & \widehat{\ST}(\D^\XX_A) \ar[d]^{\rho_\XX} \\
\widehat{W}_A \ar@{=}[rr] & & \widehat{W}_A 
}
\end{equation*}
where $\rho_\XX$ is defined by 
\begin{equation}\label{eq : extended Braid group and extended affine Weyl group}
\rho_\XX \colon \widehat{\ST}(\D^\XX_A) \longrightarrow \widehat{W}_A, \quad 
\begin{cases}
\T_{E^\XX_v} \mapsto r_{\alpha_v}, \\
\T^\L_S \mapsto (t_\delta, \gamma).
\end{cases}
\end{equation}

By \cite[Proposition 6.3]{IQ}, for a positive integer $N \ge 2$, the triangulated category $\D^N_A = \D^\XX_A / [\XX - N]$ is $N$-Calabi--Yau.
The group $\widehat{\ST}(\D^N_A)$ is also naturally defined as a subgroup of $\Aut(\D^N_A)$.
It was conjectured that the map $\iota_\XX \colon \widetilde{G}_A \longrightarrow \ST^\circ(\D^\XX_A)$ is an isomorphism.
By the commutative diagram \eqref{eq : Artin and braid}, this conjecture is equivalent to the following one.
\begin{conj}\label{conj : Artin and braid}
The map $\widehat{\iota}_\XX \colon \widehat{G}_A \longrightarrow \widehat{\ST}(\D^\XX_A)$ is an isomorphism.
\end{conj}
By \cite{IUU, Q2, W} (cf.~\cite[Corollary 6.4]{IQ}), Conjecture \ref{conj : Artin and braid} holds for $A = (1, p, q)$ and a positive integer $N \ge 2$.

Bridgeland \cite{Br} introduced the notion of a stability condition on a triangulated category.
For a given triangulated category $\D$, a stability condition is a pair $(Z, \P)$ consisting of a group homomorphism $Z \colon K_0(\D) \longrightarrow \CC$ and a $\RR$-graded family $\P$ of additive subcategories in $\D$.
He proved that the space of stability conditions on $\D$, denoted by $\Stab(\D)$, has a complex structure, which is given by the map $\Z \colon \Stab(\D) \longrightarrow \Hom_\ZZ (K_0(\D), \CC)$ defined by $(Z,\P) \mapsto Z$.
As a generalization, Ikeda--Qiu introduced $q$-stability conditions on an $\XX$-Calabi--Yau category $\D_\XX$ and considered the set $\mathrm{QStab}_s(\D_\XX)$ of $q$-stability conditions with respect to $s \in \CC$, which is also a complex manifold.

From the viewpoint of mirror symmetry, it is expected that $\Stab(\D_A)$ has the Frobenius structure given in Corollary \ref{cor : Frobenius structure} (cf.~\cite[Section 1.2]{IQ}).
On the other hand, it is also expected that $\mathrm{QStab}_s(\D^\XX_A)$ is the universal covering space of $X_A^\mathrm{reg}$ (cf.~\cite[Conjecture 10.14]{IQ2}). 
Based on Conjecture \ref{conj : Artin and braid}, we conjecture the following:
\begin{conj}\label{conj : stability condition}
Let $A = (a_1, a_2, a_3)$ with $\chi_A > 0$.
\begin{enumerate}
\item There is a biholomorphism
\begin{equation*}
\X_A \big/ \widehat{W}_A \cong \Stab(\D_A) \big/ \ST_1(\D_A).
\end{equation*}
\item For $s \in \CC$ with $\mathrm{Re}(s) \ge 2$, there is a biholomorphism 
\begin{equation*}
\X_A^\mathrm{reg} \big/ \widehat{W}_A \cong \mathrm{QStab}^\circ_s(\D^\XX_A) \big/ \widehat{\ST}(\D^\XX_A),
\end{equation*}
where $\mathrm{QStab}^\circ_s(\D^\XX_A)$ denotes the principal component of $\mathrm{QStab}_s(\D^\XX_A)$.
\end{enumerate}
\end{conj}
For $A = (1, p, q)$, Conjecture \ref{conj : stability condition} (1) follows from \cite[Theorem 6.2]{HKK}, \cite[Theorem 3.1]{DZ} and Theorem \ref{thm : main 1}. 
If $s \in \ZZ_{\ge 3}$ Conjecture \ref{conj : stability condition} (2) also follows from \cite[Theorem 5.6]{W}, \cite[Theorem 3.1]{DZ} and Theorem \ref{thm : main 1}.
We summarize these results as a proposition.
\begin{prop}\label{prop: affine A case}
For $A = (1, p, q)$ and $s \in \ZZ_{\ge 3}$, Conjecture \ref{conj : stability condition} holds.
\qed
\end{prop}


\end{document}